\documentclass[12pt]{article}

\usepackage{amsfonts}

\begin{document}

\renewcommand{\baselinestretch}{1.1} 
\normalsize

\newtheorem{lemma}{Lemma}[section]
\newtheorem{thm}[lemma]{Theorem}
\newtheorem{cor}[lemma]{Corollary}
\newtheorem{voorb}[lemma]{Example}
\newtheorem{rem}[lemma]{Remark}
\newtheorem{prop}[lemma]{Proposition}

\newfont{\gothic}{eufm10 scaled\magstep1}
\newfont{\gothics}{eufm7 scaled\magstep1}
\newfont{\gothicss}{eufm5 scaled\magstep1}

\newcommand{\gota}{\mbox{\gothic a}}
\newcommand{\gotc}{\mbox{\gothic c}}
\newcommand{\gote}{\mbox{\gothic e}}
\newcommand{\gotg}{{\mbox{\gothic g}}}
\newcommand{\gothh}{\mbox{\gothic h}}
\newcommand{\goti}{\mbox{\gothic i}}
\newcommand{\gotk}{\mbox{\gothic k}}
\newcommand{\gotm}{\mbox{\gothic m}}
\newcommand{\gotn}{\mbox{\gothic n}}
\newcommand{\gotp}{\mbox{\gothic p}}
\newcommand{\gotq}{\mbox{\gothic q}}
\newcommand{\gotr}{\mbox{\gothic r}}
\newcommand{\gots}{\mbox{\gothic s}}
\newcommand{\gott}{\mbox{\gothic t}}
\newcommand{\gotu}{\mbox{\gothic u}}
\newcommand{\gotv}{\mbox{\gothic v}}
\newcommand{\gotB}{\mbox{\gothic B}}
\newcommand{\gotG}{\mbox{\gothic G}}
\newcommand{\gotT}{\mbox{\gothic T}}
\newcommand{\gotw}{\mbox{\gothic w}}

\newcommand{\sgotg}{{\mbox{\gothics g}}}
\newcommand{\sgoth}{{\mbox{\gothics h}}}
\newcommand{\sgoti}{{\mbox{\gothics i}}}
\newcommand{\sgotm}{{\mbox{\gothics m}}}
\newcommand{\sgotn}{{\mbox{\gothics n}}}
\newcommand{\sgotr}{{\mbox{\gothics r}}}
\newcommand{\sgots}{{\mbox{\gothics s}}}
\newcommand{\sgott}{{\mbox{\gothics t}}}
\newcommand{\sgotq}{{\mbox{\gothics q}}}
\newcommand{\sgotw}{{\mbox{\gothics w}}} 
\newcommand{\ssgotm}{{\mbox{\gothicss m}}}
\newcommand{\ssgotn}{{\mbox{\gothicss n}}}

\newcommand{\Ni}{{\mathbb{N}}}
\newcommand{\Ri}{{\mathbb{R}}}
\newcommand{\Ci}{{\mathbb{C}}}
\newcommand{\Ti}{{\mathbb{T}}}
\newcommand{\Zi}{{\mathbb{Z}}}
\newcommand{\Fi}{{\mathbb{F}}}
\newcommand{\Ai}{{\mathbb{A}}}
\newcommand{\Di}{{\mathbb{D}}}

\newcommand{\proof}{\mbox{\bf Proof.} \hspace{3pt}}
\newcommand{\remark}{\mbox{\bf Remark} \hspace{3pt}}

\newcommand{\ruimte}{\vskip10.0pt plus 4.0pt minus 6.0pt}

\newcommand{\ad}{{\mathop{\rm ad}}}
\newcommand{\Ad}{\mathop{\rm Ad}}
\newcommand{\arccot}{\mathop{\rm arccot}}
\newcommand{\codim}{\mathop{\rm codim}}
\newcommand{\RRe}{\mathop{\rm Re}}
\newcommand{\IIm}{\mathop{\rm Im}}
\newcommand{\HS}{{\mathop{\rm HS}}}
\newcommand{\Tr}{{\mathop{\rm Tr}}}
\newcommand{\supp}{\mathop{\rm supp}}
\newcommand{\sgn}{\mathop{\rm sgn}}
\newcommand{\esssup}{\mathop{\rm ess\,sup}}
\newcommand{\essinf}{\mathop{\rm ess\,inf}}
\newcommand{\liminferior}{\mathop{\rm lim\, inf}}
\newcommand{\Leibniz}{\mathop{\rm Leibniz}}
\newcommand{\lcm}{\mathop{\rm lcm}}
\newcommand{\mod}{\mathop{\rm mod}}
\newcommand{\spann}{\mathop{\rm span}}
\newcommand{\ubar}{\underline{\;}}

\newcommand{\cf}{{\cal F}}
\newcommand{\ch}{{\cal H}}
\newcommand{\ck}{{\cal K}}
\newcommand{\cl}{{\cal L}}
\newcommand{\cm}{{\cal M}}
\newcommand{\co}{{\cal O}}
\newcommand{\cs}{{\cal S}}
\newcommand{\cx}{{\cal X}}
\newcommand{\cy}{{\cal Y}}
\newcommand{\cb}{{\cal B}}
\newcommand{\cu}{{\cal U}} 
\newcommand{\ct}{{\cal T}} 

\newfont{\fontcmrten}{cmr10}
\newcommand{\slbrl}{\mbox{\fontcmrten (}}
\newcommand{\slbrr}{\mbox{\fontcmrten )}}

\newfont{\fontcmsyseventeen}{cmsy10 scaled\magstep3}
\newcommand{\biggertimes}{{\mbox{{\fontcmsyseventeen \symbol{'002}}}}}

\newenvironment{remarkn}{\begin{rem} \rm}{\end{rem}}
\newenvironment{exam}{\begin{voorb} \rm}{\end{voorb}}

\begin{center}
{\bf PROPERTIES OF CENTERED RANDOM WALKS} \\[2mm] 
{\bf ON LOCALLY COMPACT GROUPS AND LIE GROUPS} \\[5mm]
Nick Dungey  \\[3mm]
Department of Mathematics, Macquarie University, NSW 2109, Australia  \\
E-mail: ndungey@ics.mq.edu.au
\end{center}

\bigskip

\noindent ABSTRACT.
The basic aim of this paper is to study asymptotic properties 
of the convolution powers 
$K^{(n)} = K*K* \cdots *K$ of a possibly non-symmetric probability density 
$K$ on a locally compact, compactly generated group $G$.  
If $K$ is centered, 
we show that the Markov operator $T$ associated with $K$ 
is analytic in $L^p(G)$ for $1<p<\infty$, 
and establish Davies-Gaffney estimates in $L^2$ 
for the iterated operators 
$T^n$.  
These results enable us to obtain various Gaussian bounds 
on $K^{(n)}$.  
In particular, 
when $G$ is a Lie group we 
recover and extend some estimates of Alexopoulos and of Varopoulos 
for convolution powers of centered densities 
and for the heat kernels of centered 
sublaplacians.  
Finally, 
in case $G$ is amenable,   
we discover that the properties of analyticity or 
Davies-Gaffney estimates hold only if $K$ is centered.

\bigskip 


\bigskip 


\section{Introduction and statement of results}   \label{s1}

The general aim of this paper is to study some properties of 
non-symmetric random walks on a locally compact group $G$.   
We consider a probability density $K$ on $G$
and the random walk on $G$ governed by $K$. 
Under certain assumptions on $K$, 
we prove time regularity estimates and Gaussian off-diagonal 
estimates for this walk. 
While such estimates are well known for symmetric random 
walks on groups or for symmetric continuous 
time diffusions on Lie groups 
(see, for example,  
\cite{HebSal,VSC,Blu1} and references therein),  
our results make no symmetry assumption.  

It is remarkable that our results apply on arbitrary 
locally compact, compactly generated groups without 
any special assumptions on the group structure 
(though it turns out that the results are most interesting 
for amenable groups).
However, in much of the paper we deal with the case where 
$G$ is a connected Lie group, since the arguments in this case 
are technically easier and of independent interest.  
The extension to arbitrary compactly generated groups 
is accomplished in the last part of the paper.   

In the Lie group case,  
we recover and extend 
some results of Alexopoulos \cite{Ale3,Ale6,Ale4} and 
Varoupoulos \cite{Var30} for non-symmetric random walks and heat kernels 
on Lie groups. 
(See also \cite{Mel} for the special case of nilpotent Lie groups.)   
Note that the analyses of \cite{Ale6,Ale4,Var30} depend on 
probabilistic methods, while our approach is completely 
different and is functional-analytic rather than probabilistic. 
Indeed, our approach essentially generalizes \cite{Dun20} 
where the case of a discrete group $G$ was studied.  
The advantage of our functional-analytic approach, 
in comparison to 
methods of \cite{Ale6,Ale4,Var30}, 
is that it does not rely on detailed knowledge of the 
structure theory or geometry 
of $G$,  
so that it applies when such structure theory is not available 
(for example, on discrete groups). 

To explain our results more precisely,  
recall 
(cf. \cite{CSal,Blu1})  
that a bounded linear operator 
$S\in {\cal L}(X)$, 
where $X$ is a complex Banach space, 
is said to be {\it analytic} 
if there exists a $c>0$ 
such that 
\[
\|(I-S)S^n \| \leq c n^{-1} 
\]
for all $n\in\Ni = \{1,2,3, \ldots\}$. 
This notion is a discrete time analogue of the usual notion 
of analyticity for a continuous time semigroup 
$(e^{tA})_{t\geq 0}$ 
which corresponds to an estimate 
$\| e^{tA}\| \leq ct^{-1}$, 
$t>0$. 
We establish the fundamental result that the Markov  
operator $T$, 
associated with the density $K$ on  
$G$, 
is analytic in $L^2(G)$ 
whenever $K$ is centered
(the definition of ``centered'' is given below).  
This theorem provides a large and interesting 
class of examples of non-self-adjoint analytic operators.  

Our second fundamental theorem gives ``Davies-Gaffney estimates'', 
that is, 
$L^2$ off-diagonal estimates, 
for the iterated operators $T^n$ 
when $K$ is centered. 
This enables us to deduce Gaussian estimates for the 
$n$-th convolution powers 
$K^{(n)}$ of $K$.  
Similar Gaussian estimates  
occur
for symmetric random walks and heat kernels in, for example,   
\cite{Var10,Car,HebSal,VSC,Mus}, 
or for non-symmetric walks on Lie groups in \cite{Ale6,Ale4,Var30}.  
But for non-symmetric walks on general groups,  
our Gaussian estimates are apparently new 
and not deducible from previously known results.

When $G$ is amenable we discover that   
the above properties, 
that is, 
analyticity of $T$ or the Davies-Gaffney estimates, 
hold if and only if the density $K$ is centered.

As already mentioned, 
our results and methods for locally compact groups 
are based on those of \cite{Dun20} for discrete groups. 
The present paper is technically more difficult than \cite{Dun20},  
since the topology of $G$ leads to problems which do not 
occur in the discrete case.

To state our results precisely we fix some 
ideas and notation 
(for relevant background material, see \cite{VSC}).  
In general, 
$G$ will denote a compactly generated, locally compact 
group and   
$dg$ will be a left invariant Haar measure on 
$G$. 
Let 
$d\widehat{g}= d(g^{-1}) = \Delta(g^{-1}) dg$
denote right Haar measure on $G$, 
where 
$\Delta\colon G\to (0,\infty)$ 
is the modular function.  

Let $K$ be a probability density on $G$ with respect to $dg$, 
that is, 
$K\colon G\to [0,\infty)$ is a Borel measurable function 
with 
$\int_G dg\, K(g) =1$.  

The density $K$ determines a random walk on $G$, 
whose distribution after $n$ steps is given by the 
$n$-th convolution power
$K^{(n)}:= K*K* \cdots *K$, 
$n\in \Ni$.  
Here,  in general the convolution of two functions 
$f_1$, $f_2$ on $G$ is given by 
\[
(f_1*f_2)(g) = \int_G dh\, f_1(h) f_2(h^{-1}g) 
\]
for all $g\in G$. 
The Markov operator 
$T$ associated with $K$ is defined by 
\[
Tf: = K*f 
\]
for any 
$f\in L^p: = L^p(G;dg)$, 
$1\leq p\leq \infty$.  
Then 
$\|T\|_{p\to p} \leq 1$ 
for all 
$p\in [1,\infty]$, 
where $\|\cdot\|_{p\to q}$ 
denotes the norm of a bounded linear operator from 
$L^p$ to $L^q$.  
Note that 
$T^n f = K^{(n)}*f$ 
for all 
$n\in\Ni$.    

Let 
$(f_1,f_2): = \int_G dg\, f_1 \overline{f_2}$ 
denote the inner product of complex-valued functions in 
$L^2= L^2(G;dg)$.  
It is easy to see that $T$ is self-adjoint 
(that is, 
$(Tf_1,f_2) = (f_1, Tf_2)$, $f_1,f_2\in L^2$)  
if and only if $K$ is symmetric in the sense that 
$K(g) = \Delta(g^{-1})K(g^{-1})$, $g\in G$, 
but we shall not assume symmetry.   

Let $L=L_G$ be the left regular representation of $G$ 
so that 
$(L(g)f)(h) = f(g^{-1}h)$ for a function 
$f\colon G\to \Ci$ 
and $g,h\in G$.   
Define the difference operators 
$\partial_g: = L(g)-I$ for 
$g\in G$.   
By definition, 
the 
``discrete Laplacian'' corresponding to 
$K$
is the operator 
\begin{equation} 
I - T = I - \int_G dg\, K(g)L(g) = 
- \int_G dg\, K(g) \partial_g. 
\label{elapdef} 
\end{equation}

Since $G$ is compactly generated, 
we may define a distance on $G$ in the following standard way.     
For the rest of the paper, 
we fix a relatively compact, 
open neighborhood 
$U$ of the identity $e$ of $G$ 
which is symmetric 
(that is, 
$U=U^{-1}$) 
and generates 
$G$.  
Then 
$G=\bigcup_{n=1}^\infty U^n$
where 
$U^n: = \{ g_1 g_2 \cdots g_n \colon g_1, \ldots, g_n \in U\}$, 
and we define 
$\rho=\rho_U \colon G\to \Ni$  
by
\[
\rho(g): = \inf\{n\in\Ni\colon g\in U^n\}, 
\;\;\; 
g\in G.  
\]
One has  
$\rho(g)=\rho(g^{-1})$ 
and $\rho(gh)\leq \rho(g)+\rho(h)$ 
for all 
$g,h\in G$.  
It is well known that the behaviour of $\rho$ 
is essentially independent of the choice of 
$U$ 
(precisely, 
if 
$\rho'$ is the distance associated with 
another such neighborhood 
$U'$, 
then 
$c^{-1} \rho \leq \rho' \leq c \rho$ for some constant 
$c>1$).  

Recall that 
$G$ is said to have polynomial 
volume growth of order $D$ 
($D\geq 0$)   
if one has 
an estimate 
\[
c^{-1} n^D \leq dg(U^n) \leq c n^D
\]
for all $n\in\Ni$.
Such groups are necessarily unimodular.  
Alternatively, 
if there exists $a>0$ with 
$dg(U^n) \geq a e^{an}$, $n\in\Ni$, 
then $G$ is said to have exponential volume growth.     
Every connected Lie group is either of polynomial 
growth or exponential growth (see \cite{Gui1}), 
but there exist examples of 
finitely generated 
discrete groups whose growth is intermediate, that is, 
neither polynomial 
nor exponential 
(\cite{Grigc}).

Let $L^p_{loc}$ be the space of functions on $G$ 
which are locally in $L^p$.   
We will often make use of the 
``$L^p$ gradient'' 
$\Gamma_p(f)$ 
of a function 
$f \in L^p_{loc}(G)$,  
defined by   
\begin{equation} 
\Gamma_p(f) := \sup_{u\in U} \| \partial_u f\|_p  
\in [0, +\infty].       
\label{egraddef}
\end{equation}

To study ``off-diagonal'' properties of the random walk, 
and in particular to obtain Gaussian estimates 
on $K^{(n)}$, 
our main technical tool will be the perturbations 
$T_{\lambda}= e^{\lambda \psi} T e^{-\lambda\psi}$, 
$\lambda\in\Ri$, 
of 
the Markov operator 
$T$. 
Here, $\psi\colon G\to \Ri$ is a suitable function
and the operator 
$T_{\lambda}$ 
is formally defined by 
\[
T_{\lambda}f 
= e^{\lambda\psi} T (e^{-\lambda\psi}f),     
\;\;\; 
f\in L^p. 
\]
The 
$\{T_{\lambda}\}_{\lambda\in\Ri}$ may be called 
``Davies perturbations'' of $T$, 
in analogy with the standard Davies perturbation of semigroups 
generated by differential operators 
(for the latter see \cite{Dav2,Dav10}). 

The class of $\psi$ we consider is the 
class
${\cal E}= {\cal E}_{G,U}$ 
of all locally bounded 
Borel measurable functions 
$\psi\colon G\to \Ri$ 
which satisfy 
\[
\Gamma_{\infty}(\psi) 
= \sup_{u\in U} \|\partial_u f\|_{\infty} \leq 1. 
\]
Note that 
$\psi$ 
itself is not assumed to be bounded.  
It is clear from the triangle inequality 
$\rho(gh)\leq \rho(g)+\rho(h)$, $g,h\in G$, 
that 
$\rho\in {\cal E}$. 

In Subsections~\ref{s1s1} to \ref{s1s3} below, 
we state our main results. 
The organization of the rest of the paper is described 
in Subsection~\ref{s1s5}.

\subsection{Densities on Lie groups}  \label{s1s1}  

To state our main results for Lie groups, in this subsection 
we shall assume that 
$G$ is a connected Lie group. 
In this case we  
usually make the following  
assumptions on the density $K$.  

(i) There exists a neighborhood $W$ of the identity $e$ 
of $G$ such that 
$\inf \{K(g)\colon g\in W\}>0$.  

(ii) For every $k>0$, one has 
\[
\int_G dg\, K(g) e^{k\rho(g)} < \infty.
\]
Assumption~(ii) requires that $K$ 
decreases sufficiently fast at infinity.   
For example, 
let us say that 
$K$ is a {\it sub-Gaussian} density 
if $K$ satisfies Assumption~(i) 
and an upper Gaussian estimate of 
the form 
\begin{equation} 
K(g) \leq c e^{-b \rho(g)^2}  
\label{ekgauss1}
\end{equation} 
for some constants 
$c,b>0$ and all 
$g\in G$. 
Then a sub-Gaussian density satisfies Assumption~(ii), 
by the well known fact that 
$\int_G dg\, e^{-\delta \rho(g)^2}$ 
is finite for any 
$\delta>0$.  
Note that any compactly supported, bounded density 
which satisfies (i) is sub-Gaussian.

The concept of centeredness for densities on a Lie group  
is defined as follows 
(cf. \cite{Ale6,Var30}).    
Let $G_0$ be the closure in 
$G$ of the commutator subgroup 
$[G,G]$,   
and consider the canonical homomorphism 
$\pi_0\colon G\to G/G_0$.   
Observe  that 
$G/G_0$ is a connected abelian Lie group and 
therefore it can be written in the form 
\[
G/G_0 \cong \Ri^{q} \times \Ti^r 
\]
for some 
$q,r\in \{0,1,2, \ldots\}$, 
where 
$\Ti:= \Ri/\Zi$. 
Define $G_1: = \pi_0^{-1} (\{0\}\times \Ti^r)$.   
Then $G_1$ is a closed normal subgroup of $G$,  
with $[G,G]\subseteq G_1$  
and 
$G_1/G_0 \cong \Ti^r$, 
$G/G_1 \cong \Ri^q$.   
Consider the homomorphism 
$\pi_1\colon G\to G/G_1 \cong \Ri^q$ 
and let 
$\pi_1^{(i)}\colon G\to \Ri$ 
be the $i$-th component of 
$\pi_1$, 
for $i\in \{1, \ldots, q\}$. 
We say that the density $K$ is 
{\it centered}  if the first order moments of 
the projection of $K$ onto $\Ri^q$ vanish, 
that is, if  
\begin{equation} 
\int_G dg\, K(g) \pi_1^{(i)}(g) = 0 
\label{edefcent} 
\end{equation} 
for all
$i\in \{1, \ldots, q\}$.

It is not difficult to show that 
$K$ is centered if and only if 
$K * \eta = \eta$ 
for all continuous homomorphisms 
$\eta\colon G\to \Ri$ 
(observe that any such homomorphism vanishes on 
$G_1$,  
and must therefore be a linear combination of the 
homomorphisms 
$\pi_1^{(i)}$.)   

The following lemma 
(whose proof is just as in \cite[Section~1.11]{Ale6}) 
shows that a general density 
is conjugate,  via a multiplicative function,   
to a centered density.

\begin{lemma}  \label{lconjugate} 
Let $K$ be a non-centered density on $G$ satisfying 
Assumptions~(i) and (ii) above.  
There exist a centered density $K'$, 
a smooth multiplicative function 
$\chi\colon G\to (0,\infty)$ 
(that is, 
$\chi(gh) = \chi(g)\chi(h)$ 
for all 
$g,h\in G$) 
and a constant 
$\delta\in (0,1)$ 
such that 
\[
K(g) = \delta \chi(g) K'(g) 
\]
for all $g\in G$. 
Setting $T'f: = K'*f$, 
we have the relations 
$K^{(n)}(g) = \delta^n \chi(g) (K')^{(n)}(g)$ 
and 
$T^n f = \delta^n \chi (T')^n (\chi^{-1}f)$ 
for all $n\in\Ni$, $g\in G$ 
and 
$f\in L^p$, $1\leq p\leq \infty$.  
\end{lemma} 

In this sense, the study of general densities reduces to 
the study of centered densities, 
and in what follows we concentrate mainly on the centered case.

We now state our first main result.   

\begin{thm} \label{tcentanal1}  
Let $K$ be a centered density on $G$ satisfying 
Assumptions~(i) and (ii) above.  
Then 
$T$ is analytic in 
$L^2$, that is, 
there exists $c>0$ with 
$\|(I-T)T^n\|_{2\to 2} \leq c n^{-1}$     
for all $n\in\Ni$.
\end{thm}

The following result is crucial for the proof of 
Theorem~\ref{tcentanal1}.

\begin{thm} \label{tgradupp1} 
Let $K$ be a centered density on $G$ satisfying Assumptions~(i) 
and (ii).  
Then there exists a  
$c>0$ 
such that 
$|((I-T)f,f)| \leq c \Gamma_2(f)^2$  
for all 
$f\in L^2$.   
\end{thm}

A theorem of Blunck \cite{Blu2} 
states that if an operator 
$S\in {\cal L}(L^{p_1})\cap {\cal L}(L^{p_2})$ 
is analytic in $L^{p_1}$ and power-bounded 
in both $L^{p_1}$ and $L^{p_2}$ 
(power-bounded means that 
$\sup_{n\in\Ni} \|S^n\| <\infty$),    
then $S$ is analytic in $L^q$ for any $q$ strictly between 
$p_1$ and $p_2$.  
Since 
$\|T^n\|_{p\to p} \leq 1$ for all 
$p\in [1,\infty]$, 
this theorem implies the following corollary of 
Theorem~\ref{tcentanal1}.

\begin{cor} \label{ccentanalp} 
Under the hypotheses of Theorem~\ref{tcentanal1}, 
then 
$T$ 
is analytic in 
$L^p$ for each 
$p\in (1,\infty)$.  
\end{cor}

Our next result is the fundamental $L^2$ off-diagonal estimate 
for centered sub-Gaussian densities.

\begin{thm} \label{tcentgauss1} 
Let $K$ be a sub-Gaussian density on $G$ 
which is centered.
There exists a $\omega>0$ such that, 
for all 
$\psi\in {\cal E}$,   
setting 
$T_{\lambda} = e^{\lambda\psi} T e^{-\lambda\psi}$ 
we have 
\[
\|T_{\lambda}^n\|_{2\to 2} \leq e^{\omega \lambda^2 n} 
\]
for all 
$n\in\Ni$ and 
$\lambda\in\Ri$.   
\end{thm} 

From Theorem~\ref{tcentgauss1} 
one can deduce the following estimates.   
Define the distance between two subsets 
$E,F$ of $G$  
by 
$d(E,F) = \inf\{\rho(gh^{-1})\colon g\in E, h\in F\}$. 
Let 
$\chi_E$ denote the characteristic function of 
$E\subseteq G$ 
(thus $\chi_E(g)=1$ or $0$ according as $g\in E$ or $g\notin E$), 
and also denote by $\chi_E$ the operator of pointwise multiplication 
$f\mapsto \chi_E f$.

\begin{thm} \label{tcentgset} 
Let $K$ be a sub-Gaussian density which is centered.  
Given any $p\in (1,\infty)$, 
there exist positive constants $c=c(p)$, $b=b(p)$ depending on $p$ 
such that 
\begin{equation} 
\| \chi_E T^n \chi_F \|_{p\to p} \leq c e^{-b d(E,F)^2/n} 
\label{egsetp} 
\end{equation} 
for all 
$n\in\Ni$ and all non-empty Borel measurable sets 
$E,F\subseteq G$.  
\end{thm} 

We refer to the estimates of Theorems~\ref{tcentgauss1} 
and 
\ref{tcentgset} as $L^2$ (or $L^p$) off-diagonal estimates 
or Davies-Gaffney estimates.  
(The cases $p=1$ or $p=\infty$ are discussed in the Remarks 
at the end of this subsection.)   
Analogous estimates are well known for symmetric random 
walks 
(see for example \cite{HebSal,CGZ}) 
or for heat semigroups on manifolds 
(for example, \cite{Dav10,Grig}).  

The next theorem shows that  
when $G$ is amenable, 
the $L^2$ 
estimates of Theorems~\ref{tcentanal1}, 
\ref{tcentgauss1} or \ref{tcentgset} 
characterize centered densities.  
(An analogous result on discrete groups  
was established in \cite{Dun20}.)   

\begin{thm} \label{tamen} 
Suppose $G$ is amenable and $K$ is a sub-Gaussian density. 
Then the following conditions~(I) to (VI) are equivalent. 

\noindent {\rm (I)}$\;$  $K$ is centered. 

\noindent {\rm (II)}$\;$  $T$ is analytic in $L^2$. 

\noindent {\rm (III)}$\;$ There exists $c>0$ such that 
\[
|((I-T)f,f)| \leq c \Gamma_2(f)^2  
\]
for all 
$f\in L^2$.  
 
\noindent {\rm (IV)}$\;$  
For each 
$\psi\in {\cal E}$, 
there exist $c,\omega>0$ such that setting 
$T_{\lambda} = e^{\lambda\psi} T e^{-\lambda\psi}$ 
we have 
\[
\|T_{\lambda}^n \|_{2\to 2} \leq c e^{\omega\lambda^2 n} 
\]
for all $n\in\Ni$ and $\lambda\in [-1,1]$.

\noindent {\rm (V)}$\;$ 
For each 
$\psi\in {\cal E}$, 
setting 
$T_{\lambda} = e^{\lambda\psi} T e^{-\lambda\psi}$ 
one has 
\[
\|T_{\lambda}\|_{2\to 2} = 1 + O(\lambda^2)
\]
for all 
$\lambda\in [-1,1]$.  

\noindent {\rm (VI)}$\;$ 
There exist $c,b>0$ such that 
\[
\| \chi_E T^n \chi_F \|_{2\to 2} \leq c e^{-b d(E,F)^2/n} 
\]
for all $n\in\Ni$ and non-empty Borel measurable 
sets $E,F\subseteq G$. 
\end{thm} 

In stating Theorem~\ref{tamen}  
we used the standard notation 
$O(\lambda^k)$, 
$\lambda\in J$, 
to denote a function 
$s=s(\lambda)$ 
satisfying an estimate 
$|s(\lambda)| \leq c |\lambda^k|$ 
for all $\lambda$ in an interval 
$J\subseteq \Ri$. 

\bigskip 

\noindent {\bf Remark.}  
It is well known that,  
for a density $K$ on $G$ satisfying 
Assumption~(i), 
one has $\|T\|_{2\to 2}<1$ 
if and only if $G$ is not amenable 
(see, for example, \cite[pp.140-142]{Pat}).   
When $G$ is not amenable, 
using $\|T\|_{2\to 2}<1$ 
it is rather easy to see that Conditions~(II) through (VI) 
of Theorem~\ref{tamen} hold even for non-centered 
densities. 
(For example, Conditions~(IV) and (V) follow from 
$\|T\|_{2\to 2}<1$ 
and the perturbation estimate of Lemma~\ref{lcrudest} below.)  
Therefore,    
our results in $L^2$ are mainly of interest for amenable 
groups.     
 
\bigskip

We proceed to state a number of results which are essentially derivable from  
the fundamental
Theorems~\ref{tcentanal1} and  \ref{tcentgauss1}.   
The following result was proved  
on amenable Lie groups
in \cite{Var30},  
by quite different methods 
(and for a slightly smaller class of densities 
$K$).

\begin{thm} \label{tcentgauss2} 
(\cite{Var30}.)    
Let $K$ be a sub-Gaussian density which is centered.  
There exist $c,b>0$ such that 
\[
K^{(n)}(g) \leq c \Delta(g)^{-1/2} 
e^{-b\rho(g)^2/n} 
\]
for all $n\in\Ni$ and $g\in G$.  
\end{thm}

We have the following analyticity and spatial regularity estimates 
for the perturbed operators 
$T_{\lambda}$.

\begin{thm} \label{tcentdiff1} 
Let $K$ be sub-Gaussian and centered. 
Set 
$T_{\lambda}:= e^{\lambda\psi} T e^{-\lambda\psi}$  
where 
$\lambda\in\Ri$ 
and 
$\psi\in {\cal E}$.  
There exist constants  
$c,\omega>0$ independent of 
$\psi$,  
such that
\[
\| (I-T_{\lambda})T_{\lambda}^n \|_{2\to 2} 
= \| e^{\lambda\psi} (I-T)T^n e^{-\lambda\psi} \|_{2\to 2} 
\leq cn^{-1} e^{\omega \lambda^2 n},  
\]
\[
\| e^{\lambda\psi} \partial_h T^n e^{-\lambda\psi}\|_{2\to 2} 
 + \| \partial_h T_{\lambda}^n \|_{2\to 2} 
\leq c\rho(h) n^{-1/2} e^{\omega \lambda^2 n},  
\]
for all $n\in\Ni$, $\lambda\in\Ri$, 
and $h\in G$ satisfying 
$\rho(h) \leq n^{1/2}$.  
\end{thm} 

In case $G$ has polynomial growth, 
Theorems~\ref{tcentgauss1} and \ref{tcentdiff1}, 
together with the arguments of \cite{Dun10},  
lead to a new proof of the following precise 
Gaussian estimates. 
These estimates were established in \cite{Ale6} 
for compactly supported densities, 
by quite different arguments involving homogenization theory.

\begin{thm} \label{tpolygauss} 
Suppose $G$ is a Lie group of polynomial volume growth of order $D$.  
Let $K$ be sub-Gaussian and centered.   
Then there are $c,b>0$ such that 
\[
K^{(n)}(g) 
+ n^{1/2} |(\partial_h K^{(n)})(g)| 
+ n |K^{(n)}(g)-K^{(n+1)}(g)|  
\leq c n^{-D/2} e^{-b\rho(g)^2/n} 
\]
for all 
$n\in\Ni$, 
$g\in G$ 
and 
$h\in U$.  
\end{thm}

\noindent {\bf Remarks.}

(a)  If the Lie group $G$ is amenable and if $K$ is centered, then 
the 
$L^p$ 
off-diagonal estimate~(\ref{egsetp}) also holds in 
the cases $p=1$ and $p=\infty$.   
This fact is essentially contained in the results of \cite{Var30}. 

Indeed, the probabilistic estimate of 
\cite[inequality~(0.9)]{Var30} on an amenable Lie group
implies that   
\begin{equation} 
\int_{\rho(g)\geq r} dg\, K^{(n)}(g) \leq c e^{-b r^2/n}
\label{eescape} 
\end{equation} 
for all 
$n\in\Ni$ and 
$r\geq 1$.  
(A version of (\ref{eescape}) for symmetric heat kernels 
was also obtained in \cite{Heb1}.)      
It is straightforward to show that (\ref{eescape}) is 
equivalent to  
the case $p=1$ (or 
$p=\infty$) 
of (\ref{egsetp}). 

(b) As already mentioned in \cite{Var30}, estimate~(\ref{eescape}) 
fails on non-amenable groups,  and also fails on 
certain amenable {\it discrete} groups of exponential growth.  
Consequently, 
the case $p=1$ of (\ref{egsetp}) fails on such groups.    

This means that 
the $L^1$ case  of (\ref{egsetp}) is essentially stronger 
than the $L^2$ case,  
for we shall see 
(see Subsection~1.4) 
that the $L^2$ case holds 
on arbitrary compactly generated groups.  

Note finally that (\ref{eescape}) holds  
when 
$G$ is a Lie group of polynomial growth, 
by a standard integration 
of the Gaussian estimate of Theorem~\ref{tpolygauss}.

(c) When $G$ has polynomial volume growth, 
the hypothesis in Theorem~\ref{tcentanal1} that 
Assumption~(ii) holds can be replaced by a weaker, 
polynomial decay condition on $K$. 
For details, see
the Remarks at the end of Section~\ref{s3}.

\subsection{Time inhomogeneous walks on Lie groups}  \label{s1s15} 

Let $G$ be a connected Lie group.   
In this subsection, we observe that some of the above results  
extend to  time inhomogeneous 
random walks on $G$.  
That is, 
given some sequence 
$(K_n)_{n=1}^{\infty}$ of densities on $G$,  
one may consider a random walk whose distribution 
after $n$ steps is given by the convolution  
\[
K_n * K_{n-1} * \cdots * K_1,  
\;\;\; n\in\Ni. 
\]
We have the following estimates 
(which generalize Theorems~\ref{tcentgauss1} and \ref{tcentgauss2})   
provided that the 
$K_n$ satisfy uniform assumptions.  

\begin{thm} \label{tunifgauss} 
Let 
$(K_n)_{n=1}^{\infty}$  
be a sequence of centered densities on 
$G$. 
Assume that  
(i) there exists a neighborhood $W$ of $e$ 
with 
$\inf\{K_n(g)\colon g\in W, n\in\Ni\}>0$,  
and 
(ii) there are $c,b>0$ with 
$K_n(g) \leq c e^{-b\rho(g)^2}$ 
for all 
$g\in G$, 
$n\in\Ni$.   
Define  
$T_n f: = K_n *f$, 
$f\in L^p$.     
Then there exist constants 
$\omega, c',b'>0$ 
such that 
\begin{equation} 
\| e^{\lambda \psi} 
T_{m+n} \cdots T_{m+1} e^{-\lambda \psi} 
\|_{2\to 2} \leq e^{\omega \lambda^2 n} 
\label{eunifg1} 
\end{equation} 
for all 
$m\in\Ni_0 =\{0,1,2, \ldots\}$, 
$n\in\Ni$, $\lambda\in\Ri$, $\psi\in {\cal E}$, 
and 
\begin{equation} 
(K_{m+n} * \cdots * K_{m+1})(g) \leq c' \Delta(g)^{-1/2} 
e^{-b'\rho(g)^2/n} 
\label{eunifg2} 
\end{equation} 
for all 
$m\in\Ni_0$, 
$n\in\Ni$ and $g\in G$. 
\end{thm}

The first estimate of Theorem~\ref{tpolygauss} extends 
to time inhomogeneous walks as follows.  

\begin{thm} \label{tpolyunif} 
Suppose that the densities 
$(K_n)_{n=1}^{\infty}$ 
satisfy the hypotheses of Theorem~\ref{tunifgauss}, 
and that $G$ has polynomial volume growth of order $D$.   
Then 
\[
(K_{m+n} * \cdots * K_{m+1})(g) \leq c n^{-D/2} e^{-b \rho(g)^2/n} 
\]
for all 
$m\in\Ni_0$ and 
$n\in\Ni$.
\end{thm}

In contrast to Theorem~\ref{tpolygauss}, 
Theorem~\ref{tpolyunif}
does not appear to be obtainable from the 
homogenization methods of \cite{Ale6}. 

Given Theorem~\ref{tpolyunif},  
then the arguments 
of \cite{Var30} probably allow one to extend 
(\ref{eescape}) to time inhomogeneous walks 
on amenable Lie groups. 
We will not, however, give the details.

In analogy with Theorem~\ref{tcentanal1},
one can ask whether estimates of type 
\[
\| (I-T_{n+1}) T_n T_{n-1} \cdots  T_1 \|_{2\to 2} 
\leq c n^{-1}, 
\;\;\; 
n\in\Ni,  
\]
are valid, 
where 
$T_j$ 
denotes the Markov operator associated 
with $K_j$. 
But the proof of Theorem~\ref{tcentanal1} 
does not yield such  estimates and we do not know when they 
are true.   

We do, though, 
have the following statement 
that a family of centered densities satisfying uniform assumptions 
is uniformly analytic.  
This result will be usefully applied to study the semigroup generated 
by a sublaplacian 
(see Subsection~\ref{s1s2} below).

\begin{thm} \label{tunifanal} 
Let 
$(K_{(\alpha)})_{\alpha\in A}$ 
be a family of centered densities on 
$G$, 
indexed by a set $A$.  
Assume that  
(i) there exists a neighborhood $W$ of $e$ 
with 
$\inf\{K_{(\alpha)}(g)\colon g\in W, \alpha\in A\}>0$,  
and 
(ii) 
$\sup_{\alpha\in A} \int_G dg\, K_{(\alpha)}(g) e^{k\rho(g)}
< \infty$ 
for each $k>0$.   
Set
$T_{(\alpha)} f = K_{(\alpha)} * f$,  
$f\in L^p$.  
Then for each 
$p\in (1,\infty)$,  
there exists a $c=c(p)>0$ 
with 
\[
\| (I-T_{(\alpha)}) T_{(\alpha)}^n \|_{p\to p} 
\leq cn^{-1} 
\]
for all $n\in\Ni$ and $\alpha\in A$.  
\end{thm}

\subsection{Sublaplacians with drift} \label{s1s2} 

In this subsection, 
$G$ will again denote a connected Lie 
group. 
An important application of the theory of Subsection~\ref{s1s1} 
above is the study of the 
semigroup and heat kernel 
associated with a sublaplacian operator on $G$.

To state some known and some new results for sublaplacians,  
let  $\gotg$ be the Lie algebra of $G$, 
consisting of all right invariant vector fields on 
$G$. 
Fix elements $A_0, A_1, \ldots, A_{d'}\in\gotg$ 
such that $A_1, \ldots, A_{d'}$ 
algebraically generate the Lie algebra $\gotg$, 
and consider the subelliptic sublaplacian 
\[
H = - \sum_{i=1}^{d'} A_i^2 + A_0 
\]
with drift term 
$A_0$. 
It is well known 
(see for example \cite[Section~IV.4]{Robm}) 
that $H$ generates a contraction 
semigroup $(e^{-tH})_{t\geq 0}$ 
in 
$L^p$, $1\leq p\leq \infty$.  
We denote by 
$K_t\colon G\to (0,\infty)$ 
the corresponding heat kernel which satisfies 
$e^{-tH} f = K_t * f$ 
and 
$K_t * K_s = K_{t+s}$  
for all 
$t,s>0$ and 
$f\in L^p$. 

For each fixed $t>0$, 
$K_t$ is a sub-Gaussian density
on $G$ 
and in fact also satisfies a lower Gaussian bound 
(see \cite[Appendix A.4]{Var6}).

One says that $H$ is {\it centered} if $A_0\in\gotg_1$,  
where $\gotg_1\subseteq \gotg$ denotes the Lie algebra 
of $G_1$. 
It is straightforward to show that 
$H$ is centered if and only if 
$K_t$ is centered for all $t>0$.   

Centered sublaplacians have been studied by various methods in 
\cite{Ale4,Var30,Mel,Dun15}.  
In this paper, 
we shall prove the following 
large time regularity result 
for the semigroup 
$e^{-tH}$,  
essentially as a consequence of Theorem~\ref{tunifanal}.

\begin{thm} \label{tsublaps} 
Let $H$ be a centered sublaplacian on $G$.  
For each 
$p\in (1,\infty)$, 
there exists 
$c=c(p)>0$ such that 
\[
\| H e^{-tH} \|_{p\to p} \leq ct^{-1} 
\]
for all $t\geq 1$. 
\end{thm} 

Using Theorems~\ref{tpolygauss} and \ref{tsublaps}, 
we will give a new proof of the following 
Gaussian estimates due to   
Alexopoulos \cite{Ale4}.

\begin{thm} \label{tsubpoly} 
(\cite{Ale4}).   
Let $G$ have polynomial growth of order $D$ and 
let $H$ be a centered sublaplacian on $G$. 
Given any right invariant vector field $X$ on $G$, 
one has estimates 
\[
K_t(g) + t^{1/2} |XK_t(g)| 
+ t \left| \frac{d}{dt}K_t(g) \right| 
\leq c t^{-D/2} e^{-b\rho(g)^2/t} 
\]
for all $t\geq 1$ and $g\in G$. 
\end{thm} 

Alternative  proofs of the first estimate 
$K_t(g) \leq c t^{-D/2} e^{-b\rho(g)^2/t}$, 
$t\geq 1$, 
of Theorem~\ref{tsubpoly} 
were given in  \cite{Dun15} and 
(in the case of nilpotent groups) 
in \cite{Mel}.   

Note that when $G$ has polynomial growth, 
a standard integration 
of the Gaussian estimate on 
$(d/dt)K_t = -H K_t$ of Theorem~\ref{tsubpoly}  
implies the estimate of Theorem~\ref{tsublaps}, 
for any 
$p\in [1,\infty]$.    
But when $G$ has exponential growth, 
the estimate of Theorem~\ref{tsublaps} appears to be new.

\subsection{Locally compact groups}   \label{s1s3} 

In this subsection,  
we allow $G$ to be any locally compact, compactly generated group. 
It is remarkable that the main results of 
Subsections~\ref{s1s1} and \ref{s1s15} extend to this 
general situation.    

For simplicity,  we make the following assumption of compact support on 
the density 
$K$.   

(i)' $K$ is compactly supported, bounded, 
and 
$\inf\{K(g) \colon g\in W\}>0$,  
where $W$ 
is some relatively compact, open and symmetric neighborhood 
of the identity $e$ of $G$ which generates $G$. 

\bigskip

The notion of centeredness is extended to this general situation
as follows.   
Set 
$G_0 = \overline{[G,G]}$ 
and consider the projection
$\pi_0\colon G\to G/G_0$. 
Since  
$G/G_0$ is a compactly generated locally compact abelian group, 
by a standard structure theorem 
(see \cite[Theorem~II.9.8]{HR}) 
it can be written as a direct product 
$G/G_0 \cong \Ri^q \times \Zi^r \times M$ 
where 
$q,r \in \{0,1,2, \ldots\}$ 
and 
$M$ is a compact abelian group.  
Set $G_1  = \pi_0^{-1}(\{0\}\times M)$. 
Then 
$G_1$ is a closed normal subgroup of $G$ 
with 
$[G,G]\subseteq G_1$, 
and 
$G_1/G_0 \cong M$ is compact.     
Consider the homomorphism 
$\pi_1\colon G\to G/G_1 \cong \Ri^q \times \Zi^r$ 
with components 
$\pi_1^{(i)}\colon G\to \Ri$, 
$i \in \{1, \ldots, q\}$, 
and 
$\pi_1^{(i)}\colon G\to \Zi$, 
$i\in \{q+1, \ldots, q+r\}$.    
We say that 
$K$ 
is {\it centered}  if  
(\ref{edefcent}) holds for each
$i\in \{1, \ldots, q+r\}$. 

Just as in the Lie group case, 
one sees that 
$K$ is centered if and only if 
$K * \eta = \eta$ 
(or,  equivalently,  
$(I-T)\eta=0$)    
for every continuous homomorphism 
$\eta\colon G\to \Ri$.

If $K$ is a non-centered density satisfying (i)', 
then there exists a centered density 
$K'$ such that 
$K = \delta \chi K'$ 
where $\delta\in (0,1)$ is a constant and $\chi\colon G\to (0,\infty)$ 
is a continuous multiplicative function. 
The proof of this fact is similar to the proof of Lemma~\ref{lconjugate} 
(cf. \cite{Ale6,Ale5}) 
and is left to the reader.  

Concerning analyticity of $T$,
we have the following extension of Theorems~\ref{tcentanal1} 
and \ref{tgradupp1}.

\begin{thm} \label{tlcanal}   
Let $G$ be a locally compact, compactly generated group and let 
$K$ be  a centered density on $G$ 
satisfying Assumption~(i)'.  
Then there exists $c>0$ such that 
\begin{equation} 
|((I-T)f_1,f_2)| \leq c \Gamma_2(f_1)\Gamma_2(f_2) 
\label{eupgd22} 
\end{equation} 
for all $f_1,f_2\in L^2$. 
Moreover, 
$T$ is analytic in $L^p$ 
for each $p\in (1,\infty)$. 
\end{thm} 

We also have off-diagonal estimates generalizing 
Theorems~\ref{tcentgauss1}, \ref{tcentgset} and 
\ref{tcentgauss2}, 
as follows.    

\begin{thm} \label{tlcgauss} 
Let $G$ be a locally compact, compactly generated group, 
and let $K$ be a centered density on 
$G$ satisfying Assumption~(i)'. 
Then the following statements hold. 

\noindent {\rm (I)}$\;$ 
There is $c>0$ such that 
\begin{equation} 
|((I-T)f_1,f_2)| \leq c \Gamma_{\infty}(f_1)\Gamma_1(f_2) 
\label{eupgdinf1} 
\end{equation} 
for all $f_1\in L^\infty$ and $f_2\in L^1$.

\noindent {\rm (II)}$\;$  
There is 
$\omega>0$ 
such that 
\[
\| e^{\lambda \psi} T^n e^{-\lambda\psi} \|_{2\to 2} 
\leq e^{\omega \lambda^2 n} 
\]
for all 
$n\in\Ni$, $\lambda\in \Ri$, $\psi\in {\cal E}$.  
For each $p\in (1,\infty)$ there exist 
$c=c(p)$, $b=b(p)>0$ 
such that 
\[
\| \chi_E T^n \chi_F \|_{p\to p} \leq c e^{-b d(E,F)^2/n} 
\]
for all $n\in\Ni$ and all Borel sets $E,F\subseteq G$, 
where 
$d(E,F)= \inf\{\rho(gh^{-1})\colon g\in E, h\in F\}$.

\noindent {\rm (III)}$\;$ 
There are $c,b>0$ such that 
\[
K^{(n)}(g) \leq c \Delta(g)^{-1/2} 
e^{-b\rho(g)^2/n} 
\]
for all $n\in\Ni$ and $g\in G$.  
\end{thm}

Moreover, 
the estimates of Theorem~\ref{tcentdiff1} hold 
when $K$ is centered.   
As a consequence, 
we can generalize  
Theorem~\ref{tpolygauss} 
for any locally compact group of polynomial growth.

\begin{thm} \label{tlcpolygauss} 
Suppose that $G$ is a locally compact, compactly generated group 
of polynomial volume growth of order $D$.  
Let $K$ be a centered density on $G$ satisfying Assumption~(i)'.  
Then there are $c,b>0$ such that 
\[
K^{(n)}(g) 
+ n^{1/2} |(\partial_h K^{(n)})(g)| 
+ n |K^{(n)}(g) - K^{(n+1)}(g)| 
\leq c n^{-D/2} e^{-b \rho(g)^2/n} 
\]
for all 
$n\in\Ni$, 
$g\in G$ 
and 
$h\in U$.  
\end{thm} 

For amenable groups 
we have the following extension of Theorem~\ref{tamen}.

\begin{thm} \label{tamengen} 
Let $G$ be an amenable locally compact, compactly generated group
and $K$ a density on $G$ satisfying Assumption~(i)'. 
Then the Conditions (I) through (VI) in Theorem~\ref{tamen} 
are all equivalent.  
\end{thm} 

If $G$ is non-amenable then the situation 
is as follows (compare the Lie group case 
in Subsection~\ref{s1s1}): 
for any,  possibly non-centered,  density $K$ satisfying Assumption~(i)',  
one has 
$\|T\|_{2\to 2}<1$ and Conditions~(II) through (VI) in Theorem~\ref{tamen} 
all hold.    

\bigskip

\noindent {\bf Remark.} 
Let us mention some further  
$L^\infty$ estimates and Gaussian estimates 
for $K^{(n)}$. 
Observe that if 
$K$ is centered and one has a uniform upper bound of the form 
\begin{equation} 
\|  \Delta^{1/2} K^{(n)} \|_{\infty} \leq \gamma(n)
\label{euppm} 
\end{equation} 
for all $n\in\Ni$ and some function
$\gamma\colon \Ni \to (0,\infty)$,  
then interpolation with the bound of Theorem~\ref{tlcgauss}, 
part~(III), yields that 
\[
\Delta(g)^{1/2}  K^{(n)}(g) 
\leq \gamma(n)^{1-\varepsilon} c^{\varepsilon} e^{-b\varepsilon 
\rho(g)^2/n}
\]
for any 
$\varepsilon\in (0,1)$, 
$g\in G$ and 
$n\in\Ni$.  
For example, 
if (\ref{euppm}) holds with 
$\gamma(n)= c_1 e^{-c_2 n^{\beta}}$ where $\beta\in (0,1]$,  
then we obtain 
$K^{(n)}(g) \leq c_3 \Delta(g)^{-1/2} e^{-c_4 n^{\beta}} e^{-c_4 \rho(g)^2/n}$ 
for some constants 
$c_3,c_4>0$.   

For example, 
if $G$ is unimodular and of exponential volume growth 
then (\ref{euppm}) holds with 
$\gamma(n)= c_1 e^{-c_2 n^{1/3}}$, 
by a well known theorem (see \cite[Chapter~VII]{VSC}). 
Therefore, in this case an estimate of type  
\[
K^{(n)}(g) \leq c e^{-bn^{1/3}} e^{-b\rho(g)^2/n} 
\]
holds for all 
$n\in\Ni$, $g\in G$,  
when $K$ is centered.   
This appears to be a new result for non-symmetric densities.

\bigskip 

Finally, 
for time inhomogeneous random walks on $G$ 
we record the following extension of 
Theorem~\ref{tlcgauss}.    

\begin{thm} \label{tlcunif} 
Let 
$G$ be a locally compact, compactly generated group 
and let 
$(K_n)_{n=1}^{\infty}$ 
be a sequence of compactly supported, centered densities on $G$. 
Define 
$T_n f: = K_n *f$, 
$f\in L^p$.   
Suppose there exist 
$\varepsilon>0$ and two relatively compact, symmetric generating 
neighborhoods $W$, $W'$ of the identity $e$ such that 
$W\subseteq W'$, 
$W$ is open, 
and 
\[
\varepsilon \chi_W \leq K_n \leq \varepsilon^{-1} \chi_{W'} 
\]
for all 
$n\in\Ni$. 
Then there exist 
$\omega,c,b>0$ 
such that 
\[
\| e^{\lambda \psi} 
T_{m+n} \cdots T_{m+1} e^{-\lambda \psi} 
\|_{2\to 2} \leq e^{\omega \lambda^2 n} 
\]
for all 
$m\in\Ni_0$, 
$n\in\Ni$, $\lambda\in\Ri$, $\psi\in {\cal E}$, 
and 
\[
(K_{m+n} * \cdots * K_{m+1})(g) \leq c \Delta(g)^{-1/2} 
e^{-b\rho(g)^2/n} 
\]
for all $m\in\Ni_0$, 
$n\in\Ni$ and $g\in G$. 
\end{thm}

\subsection{Organization of the paper}  \label{s1s5}

An outline of the paper is as follows. 

In Section~\ref{s2} we begin the proof of Theorem~\ref{tcentanal1} 
for a Lie group;  
the proof is completed in Section~\ref{s3} with 
the proof of Theorem~\ref{tgradupp1}.  
The underlying idea of the proof of Theorem~\ref{tgradupp1} 
is that, 
roughly speaking,  
if $K$ is centered then $I-T$ is a second-order difference 
operator without first-order terms, 
that is, it is a generalized linear combination 
of operators 
$\partial_{g_1} \cdots \partial_{g_k}$ 
where $g_1, \ldots, g_k\in G$ and $k\geq 2$ 
(see \cite{Dun20} for the analogous assertion on discrete groups).   
To develop this idea precisely on Lie groups, 
however,  requires some effort.

In Section~\ref{s4} we establish Theorem~\ref{tcentgauss1}.  
For this proof we need a certain variation 
of the estimate of Theorem~\ref{tgradupp1} 
(see (\ref{euppgradpsi}) below).

The ideas and results of Sections~\ref{s2}, \ref{s3} and \ref{s4} 
are fundamental for this article, and are used repeatedly 
in later sections.

Next, sections~\ref{s5} and \ref{s6} 
contain the proofs of Theorems~\ref{tcentgset} 
and \ref{tamen} respectively.  
The proof of Theorem~\ref{tcentgset} is based on 
an equivalence between two different forms of $L^p$ 
off-diagonal estimates,  and shows that 
Theorems~\ref{tcentgset} and \ref{tcentgauss1} are essentially 
equivalent.

Section~\ref{s7} describes  
the proofs of Theorems  
\ref{tcentgauss2} to \ref{tunifanal}, 
utilizing results 
in Sections~\ref{s2} to \ref{s4}.   

In Section~\ref{s8}, we consider sublaplacians on the 
Lie group $G$,  and prove 
Theorems~\ref{tsublaps} and \ref{tsubpoly}.  

Finally, in Sections~\ref{s9} and \ref{s10} 
we consider general locally compact compactly 
generated groups and derive the results stated 
in Subsection~\ref{s1s3}.  
The proofs are based on similar ideas as in the Lie group case, 
but significant difficulties occur and 
we rely on some deep structure theorems 
for locally compact groups (cf. \cite{MZ}).   

In general, 
$c, c', b$ and so on will denote positive constants 
whose value may change from line to line when convenient.

\section{Proof of Theorem~\ref{tcentanal1}}  \label{s2} 

In this section, unless otherwise stated 
$G$ will be a connected 
Lie group.   
Note, however, 
that the arguments of this section do not use Lie theory 
and extend (with only minor changes) to locally compact, 
compactly generated groups.   

To prove Theorem~\ref{tcentanal1} we will apply the following 
general characterization of analytic operators 
due to Nevanlinna 
(see \cite[Theorem~4.5.4]{N1},  
\cite[Theorem~2.1]{N2}, 
and 
\cite{Blu1,Blu2}).    
Let 
$\Di = \{\lambda \in \Ci\colon |\lambda|<1\}$ 
be the open unit disk.      

\begin{thm} \label{tanalcrit} 
Let 
$X$ be a complex Banach space and let  
$S\in {\cal L}(X)$. 
The following three conditions are equivalent.  

\noindent {\rm (I)}$\;$  
$S$ is power-bounded 
(that is, 
$\sup\{ \|S^n\|\colon n\in\Ni\} < \infty$) 
and analytic.  

\noindent {\rm (II)}$\;$ 
$(e^{-t(I-S)})_{t\geq 0}$ 
is a bounded analytic semigroup in $X$,  
and  
the spectrum of 
$S$ is contained in 
$\Di \cup \{1\}$.   

\noindent {\rm (III)}$\;$ 
There exists a 
$c>0$ such that 
$\| (\lambda I - S)^{-1} \| \leq c|\lambda -1|^{-1}$ 
for all 
$\lambda \in\Ci$ with 
$|\lambda|>1$. 
\end{thm} 

To prove Theorem~\ref{tcentanal1} 
we will verify Condition~(II) of Theorem~\ref{tanalcrit} 
when 
$S=T$ 
(where 
$Tf: = K*f$)  
and 
$X=L^2=L^2(G;dg)$.   
The following general proposition,  
and Assumption~(i),  
yield  the desired condition on the spectrum 
of $T$.  
Denote by 
$\sigma_p(S)$ the spectrum of an operator 
$S\in {\cal L}(L^p)$.

\begin{prop} \label{pspect} 
Let $G$ be a locally compact group with left Haar measure 
$dg$,  
let  
$P\colon G\to [0,\infty)$ 
be a probability density on 
$G$, 
and define the operator 
$S$ by  
$Sf = P*f$, 
$f\in L^2:=L^2(G;dg)$. 
Suppose there exists a neighborhood 
$W$ of the identity 
of $G$ such that 
$\inf\{P(g)\colon g\in W\} >0$.    
Then 
$\sigma_2(S)\subseteq \Di \cup \{1\}$. 
\end{prop} 

The proof of Proposition~\ref{pspect} uses 
the following Hilbert space result.   
For any 
$\lambda\in\Ci$, $A\subseteq \Ci$, 
let us write 
$dist(\lambda,A) := \inf\{|\lambda-a|\colon a\in A\}$.

\begin{lemma} \label{lspectr} 
For any 
$S\in {\cal L}(L^2)$, 
the spectrum  
$\sigma_2(S)$ is contained in the closure 
$\overline{\Theta(S)}$ of the set 
$\Theta(S) := \{(Sf,f)\colon f\in L^2, \|f\|_2=1\} \subseteq \Ci$.   
Moreover, 
for each 
$\lambda\in \Ci\backslash \overline{\Theta(S)}$, 
setting 
$d_{\lambda}: = dist(\lambda, \overline{\Theta(S)}) 
= dist(\lambda,\Theta(S))>0$ 
one has 
$ \|(\lambda I - S)^{-1}\|_{2\to 2} \leq d_{\lambda}^{-1}$.   
\end{lemma} 

The first statement of the lemma is standard 
(see 
\cite[Corollary~V.3.3]{Kat1}). 
The second statement then follows by observing that 
$dist(\lambda,\Theta(S)) 
\leq |\lambda-(Sf,f)| \leq \|(\lambda I-S)f\|_2$ 
for 
$\|f\|_2=1$.  

\bigskip 

\noindent {\bf Proof of Proposition~\ref{pspect}.}  
Let $P$ and 
$W$ be as in the hypothesis,   
and choose  a compact neighborhood $V$ 
of $e$ such that $V=V^{-1}$ 
and $VV\subseteq W$. 
Let 
$\Delta\colon G\to (0,\infty)$ 
be the modular function of $G$,   
so that 
$d\widehat{g}: = \Delta(g^{-1}) dg$ 
is right Haar measure on $G$. 
Set 
\[
P_1 : = \varepsilon (\Delta^{-1}\chi_V) * \chi_V, 
\]  
where 
$\chi_V$ is the characteristic function of $V$ and 
$\varepsilon>0$ is a constant chosen small enough so that 
$P_1 \leq P$.   
Define 
$S_1 f := P_1*f$   
for all 
$f\in L^2$.  
It is easy to check that 
\[
(S_1 f, f) = \varepsilon (\chi_V * f, \chi_V * f) \geq 0, 
\]
so that 
$S_1$ is non-negative self-adjoint in 
$L^2$.  
Put 
$\delta = \int_G dg\, P_1(g)$.  
Because 
$0\leq P_1\leq P$, 
we have  
$0<\delta\leq 1$  
and 
\[
\|S_1\|_{2\to 2} \leq \delta, 
\;\;\; 
\| S-S_1\|_{2\to 2} \leq \int_G dg\, (P-P_1)(g) = 1-\delta. 
\]
Let 
$\|f\|_2=1$.  
Since 
$(S_1 f,f) \in [0,\delta]$, 
then 
$(Sf,f)= (S_1 f,f) + ((S-S_1)f,f)$ 
is an element of the set 
\begin{equation} 
\Lambda_{\delta}: = 
\{\lambda \in \Ci\colon dist(\lambda , [0,\delta]) \leq 1-\delta \}.   
\label{elambdadef} 
\end{equation} 
By Lemma~\ref{lspectr}   
and since 
$\delta\in (0,1]$,  
we conclude that 
$\sigma_2(S) \subseteq \overline{\Lambda_{\delta}} 
= \Lambda_{\delta} 
\subseteq \Di\cup\{1\}$. 
\hfill $\Box$

\bigskip

By Theorem~\ref{tanalcrit} and Proposition~\ref{pspect},  
the proof of Theorem~\ref{tcentanal1} 
is reduced to showing that, when $K$ is centered,  
$e^{-t(I-T)}$ 
is a bounded analytic semigroup in 
$L^2$. 
To show this, 
by a well known semigroup result 
(see \cite[Theorem~IX.1.24]{Kat1})  
it suffices to obtain a sectorial estimate 
of form 
\begin{equation} 
|((I-T)f,f)| \leq c \RRe ((I-T)f,f)
\label{esect1} 
\end{equation} 
for all 
$f\in L^2$.  
We will begin the proof of (\ref{esect1}) in this section 
and will complete the proof in Section~\ref{s3}.

We will need the identities 
\begin{equation} 
\partial_{g_1 \cdots g_n} =  \partial_{g_1} + L(g_1) \partial_{g_2} + 
   \cdots + L(g_1 \cdots g_{n-1}) \partial_{g_n} 
\label{ediffid}  
\end{equation} 
and 
\begin{equation} 
\partial_{g_1 \cdots g_n}  =  \partial_{g_1} + \cdots + \partial_{g_n} 
+ \sum_{k=2}^n 
   \sum_{i_1, \ldots, i_k} \partial_{g_{i_1}} \ldots \partial_{g_{i_k}}
\label{ediffid2} 
\end{equation} 
for all $n\in\Ni$ and 
$g_1, \ldots, g_n\in G$.  
In the second identity,  
the inner sum is over all integers 
$i_1, \ldots, i_k$ 
with 
$1\leq i_1 <i_2 < \cdots <i_k \leq n$.     

The next lemma records basic facts about  
$\Gamma_p$ defined by (\ref{egraddef}).

\begin{lemma} \label{lgradbas} 
One has 
\begin{equation} 
\|\partial_g f\|_p  \leq  \rho(g) \Gamma_p(f) 
\label{enabla1} 
\end{equation} 
and 
\begin{equation} 
\Gamma_p (L(g)f)  \leq 3 \rho(g) \Gamma_p(f)  
\label{enabla2} 
\end{equation} 
for all $g\in G$ and $f\in L^p$. 
Given any relatively compact, Borel measurable neighborhood 
$V$ of $e$ in $G$, 
there exists a $c=c(V)>1$ such that 
\[
c^{-1} \Gamma_2(f)^2  \leq \int_V dh\, \|\partial_h f\|_2^2 
\leq c \Gamma_2(f)^2 
\]
for all $f\in L^2$. 
\end{lemma} 
\proof\
Inequality (\ref{enabla1}) follows easily by 
writing 
$g=g_1 g_2 \cdots g_n$ with 
$g_j\in U$, 
$n=\rho(g)$, 
and applying (\ref{ediffid}).

To prove (\ref{enabla2}),  
note that 
$\partial_u L(g) f 
= -\partial_g f + L(u) \partial_g f + \partial_u f$
for 
$u\in U$, $g\in G$.  
Therefore 
by (\ref{enabla1}), 
\[
\| \partial_u L(g) f\|_p 
 \leq 2 \|\partial_g f\|_p + \|\partial_u f\|_p
\leq (2\rho(g)+1)  \Gamma_p(f) \leq 3\rho(g) \Gamma_p(f) 
\]
for all $u\in U$ and $g\in G$, 
which implies (\ref{enabla2}).   

In the final statement of the lemma, 
the upper bound on 
$\int_V dh\, \|\partial_h f\|_2^2$ 
follows easily from (\ref{enabla1}). 
The lower bound follows as in the proof 
of \cite[Proposition~VII.3.2]{VSC} and 
we refer there for details. 
\hfill $\Box$ 

\bigskip 

\noindent {\bf Remark.}  
The statements of Lemma~\ref{lgradbas} actually
hold on any locally compact compactly generated 
group $G$, 
provided that 
$V$ is a relatively compact, 
open 
neighborhood of $e$ which generates $G$.     

\bigskip

The first step in the proof of (\ref{esect1}) is the following 
observation.   

\begin{lemma} \label{lrealest} 
There exists a  
$c>1$ such that 
$c^{-1} \Gamma_2(f)^2  \leq \RRe ((I-T)f,f) 
\leq c \Gamma_2(f)^2$   
for all $f\in L^2$. 
\end{lemma} 
\proof\
Let $T^*$ be the $L^2$-adjoint of $T$ and consider the 
self-adjoint operator  
$\widehat{T}: = 2^{-1}(T+ T^*)$.  
For 
$f\in L^2$, 
note that 
\[
\RRe ((I-T)f,f)= 2^{-1}((I-T)f,f) + 2^{-1}(f,(I-T)f) 
= ((I-\widehat{T})f,f)   
\]
and that 
$\widehat{T}f = \widehat{K}*f$,   
where 
$\widehat{K}$ is the probability density
\[
\widehat{K}(g) : = 2^{-1}(K(g) + \Delta(g^{-1}) K(g^{-1})), 
\;\;\; 
g\in G. 
\] 
Using 
$d(g^{-1})= \Delta(g^{-1})dg$ 
and the symmetry 
$\widehat{K}(g)=\Delta(g^{-1}) \widehat{K}(g^{-1})$, 
one checks that  
\begin{eqnarray*} 
2^{-1} \int_G dg\, \widehat{K}(g) \partial_{g^{-1}} \partial_g 
& = & 2^{-1} \int_G dg\, \widehat{K}(g)  
  (2I - L(g) - L(g^{-1})) \\
& = & I - 2^{-1} \int_G dg\, \widehat{K}(g) L(g) 
         - 2^{-1} \int_G d(g^{-1})\, \widehat{K}(g^{-1}) L(g)  \\
& = & I - \int_G dg\, \widehat{K}(g) L(g) 
= I - \widehat{T}.
\end{eqnarray*} 
Because 
$\partial_{g^{-1}}$ has adjoint 
$\partial_g$, 
then 
\begin{eqnarray*} 
((I-\widehat{T})f,f) & = & 
2^{-1} \int_G dg\, 
\widehat{K}(g) \|\partial_g f\|_2^2   \\
& \leq & 2^{-1} \int_G dg\, 
   \widehat{K}(g) \rho(g)^2  \Gamma_2(f)^2 
  \\ 
& \leq & c  \Gamma_2(f)^2  
\end{eqnarray*} 
for all $f\in L^2$, 
since 
$\int_G dg\, \widehat{K}(g) \rho(g)^2  <\infty$ 
by Assumption~(ii).  

Next, 
by Assumption~(i) there exists a compact neighborhood 
$V'$ of $e$ 
such that 
$\inf\{\widehat{K}(g)\colon g\in V'\} >0$. 
We obtain an estimate  
\[
((I-\widehat{T})f,f) =  2^{-1} \int_G dg\, 
\widehat{K}(g) \|\partial_g f\|_2^2   
\geq c \int_{V'} dg\, \|\partial_g f\|_2^2 
\geq c' \Gamma_2(f)^2,   
\]
where the last step follows by the last statement of 
Lemma~\ref{lgradbas}.  
The lemma follows.     
\hfill $\Box$

\bigskip 

\noindent {\bf Remark.} 
The proof of Lemma~\ref{lrealest}  
does not require that $K$ be centered:  
the lemma is valid for any density 
satisfying Assumptions~(i) and (ii).    

\bigskip

In Section~\ref{s3}, 
we will prove Theorem~\ref{tgradupp1}.

Combining Theorem~\ref{tgradupp1} and 
Lemma~\ref{lrealest} 
gives an estimate of form (\ref{esect1}).  
Then 
$(e^{-t(I-T)})_{t\geq 0}$ 
is a bounded analytic semigroup in 
$L^2$, 
and Theorem~\ref{tcentanal1} follows.

\section{Proof of Theorem~\ref{tgradupp1}} \label{s3}

Let $G$ be a connected Lie group with Lie algebra 
$\gotg$. 
In this section we prove Theorem~\ref{tgradupp1}, 
which will complete the proof of Theorem~\ref{tcentanal1}.  

A key element of the proof is to establish 
estimates of type  
\begin{equation} 
|(\partial_g f_1, f_2) \leq c(g)  \Gamma_2(f_1) \Gamma_2(f_2)
\label{eg1gen}
\end{equation} 
for all $f_1,f_2\in L^2$, 
where $g$ is an element of the subgroup $G_1$ of $G$ 
and 
$c(g)>0$ is a constant which may depend on $g$. 
The idea behind the proof of such estimates is that when $g\in G_1$
then $\partial_g$ is a ``second-order'' difference operator 
in the sense that it is a linear combination  
of operators $\partial_{g_1} \cdots \partial_{g_k}$ 
where $k\geq 2$, $g_1, \ldots, g_k\in G$.  
To make this idea precise, we shall apply some structure 
theory of Lie groups.  
(Later, 
by analogous but more complicated arguments we shall obtain  
(\ref{eg1gen}) on any locally compact compactly generated group; 
see Proposition~\ref{pg0g1ext}.)

Recall that 
$[G,G]$ consists of all finite products of 
the commutators 
$[g_1,g_2]: = g_1 g_2 g_1^{-1} g_2^{-1}$, 
$g_1,g_2\in G$.  
It is well known 
(see \cite[Theorem~3.18.7]{Vara})  
that 
$[G,G]$ is a connected, possibly non-closed,  Lie subgroup of 
the connected Lie group $G$, 
and has  
Lie algebra 
$[\gotg,\gotg]$. 
The groups $G_0 = \overline{[G,G]}$ and 
$G_1$,  
defined in Subsection~\ref{s1s1},  
are closed, connected Lie subgroups of $G$. 
(Connectedness of 
$G_1$ follows from the connectedness 
of 
$G_0$ and of 
$G_1/G_0 \cong \Ti^r$.)    
Denote by $\gotg_1$ the Lie algebra of $G_1$.

\begin{lemma} \label{lggnbhd} 
There exist 
$m\in\Ni$  
and a compact symmetric neighborhood 
$U'$ of $e$ in $G$ 
such that, 
setting   
\[
U'' := \{ [g_1,g_2]\colon g_1,g_2\in U'\}, 
\]
then the set 
$(U'')^m: = \{u_1 \cdots u_m \colon u_j\in U''\}$ 
is a neighborhood of the identity for the Lie group 
$[G,G]$.     
\end{lemma} 
\proof\
Let $m\in\Ni$ and consider  
the smooth mapping 
$\Psi\colon G^{2m} \to [G,G]$ 
defined by 
\[
\Psi(g_1,h_1,g_2,h_2, \ldots, g_m,h_m):  = 
[g_1,h_1] \cdots [g_m, h_m]
\]
for all 
$g_j, h_j\in G$.  
The argument of \cite[p.242]{Vara} 
shows that  
if $m$ is chosen large enough,  
then the differential  
\[
\Psi_*|_p \colon T_p(G^{2m})\to T_e([G,G])  
\]
is surjective 
at some point 
$p= (a_1,e,a_2,e, \ldots, a_m,e)\in  
G^{2m}$, 
where 
$a_1, \ldots, a_m \in G$. 
The lemma then follows by taking any  compact symmetric neighborhood $U'$ 
of $e$ 
whose interior contains $a_1, \ldots, a_m$. 
\hfill $\Box$ 

\bigskip 

\noindent {\bf Remark.}  
As an aside, 
a closer examination of the argument of \cite[pp.241-242]{Vara}
shows that one can take 
$m=\dim([\gotg,\gotg])$, 
and that 
$a_1, \ldots, a_m$ could be chosen 
arbitrarily chose to $e$; 
thus the conclusion 
of Lemma~\ref{lggnbhd} actually 
holds for any neighborhood 
$U'$ of $e$ in $G$.

\bigskip

Next, 
we need the following structural result for $G_1$.  

\begin{lemma} \label{lg1dec} 
There exists a compact,  
connected, abelian subgroup 
$C$ of $G$, 
with Lie algebra $\gotc$,   
such that 
$\gotg_1 = \gotc + [\gotg, \gotg]$ 
and 
$G_1 = C [G,G]$.
(Here, $+$ denotes a sum of vector subspaces which is not 
necessarily direct.)     

If $U',U'',m$ are as in Lemma~\ref{lggnbhd}, 
then 
$U_1 : = C (U'')^m$ 
is a neighborhood of the identity for the group $G_1$.   
\end{lemma} 
\proof\ 
The existence of $C$ with the properties 
$\gotg_1 = \gotc + [\gotg,\gotg]$, 
$G_1 = C [G,G]$,  
is proved in 
the Appendix of \cite{Dun15} 
(see also \cite[Appendix]{Var30} for similar results when $G$ 
is amenable).   
That 
$C (U'')^m$ is a neighborhood of the identity for 
$G_1$ then follows from Lemma~\ref{lggnbhd}. 
\hfill $\Box$ 

\bigskip 

In the rest of this section, 
we shall 
fix 
$U',U'',m, C$ 
and 
$U_1 = C (U'')^m$ 
with properties as in Lemmas~\ref{lggnbhd} and \ref{lg1dec}.

\begin{lemma} \label{luest} 
There exists $c>0$ such that 
$|(\partial_g f_1,f_2)| \leq 
c \Gamma_2(f_1) \Gamma_2(f_2)$ 
for all 
$g\in U''$  
and 
$f_1,f_2\in L^2$.   
\end{lemma} 
\proof\
Let $g=[g_1,g_2]\in U''$ 
where 
$g_1,g_2\in U'$.  
The identity 
\[
\partial_g = - L(g_1 g_2) 
(\partial_{g_2^{-1}} \partial_{g_1^{-1}} 
- \partial_{g_1^{-1}} \partial_{g_2^{-1}}) 
\]
yields 
\[
(\partial_g f_1,f_2) 
= - (\partial_{g_1^{-1}} f_1, \partial_{g_2} L(g_2^{-1} g_1^{-1}) f_2) 
  + (\partial_{g_2^{-1}} f_1, \partial_{g_1} L(g_2^{-1} g_1^{-1}) f_2). 
\]  
We use the Cauchy-Schwarz inequality, (\ref{enabla1}),  
and (\ref{enabla2}) to obtain 
\begin{eqnarray*} 
|(\partial_g f_1,f_2)| 
& \leq & c \Gamma_2(f_1) \Gamma_2 (L(g_2^{-1}g_1^{-1}) f_2) 
    \\
& \leq & c' \Gamma_2(f_1) \Gamma_2(f_2)
\end{eqnarray*} 
for all 
$g\in U''$ and 
$f_1,f_2\in L^2$.    
\hfill $\Box$

\begin{lemma} \label{luoneest} 
There exists 
$c>0$ such that 
$|(\partial_g f_1,f_2)| \leq 
c \Gamma_2(f_1) \Gamma_2(f_2)$  
for all 
$g\in U_1$ and  
$f_1,f_2\in L^2$.   
\end{lemma} 
\proof\
Let 
$ds$ be Haar measure on the compact group 
$C$, 
normalized so that 
$ds(C)=1$.  
For any $s_0\in C$, 
we have the operator identity 
\begin{equation} 
\partial_{s_0} = - \int_C ds\, \partial_s \partial_{s_0}.  
\label{ecpctrep} 
\end{equation} 
Using (\ref{enabla1}) we deduce that 
\[
|(\partial_{s_0}f_1,f_2)| 
\leq \int_C ds\, |(\partial_{s_0} f_1, \partial_{s^{-1}}f_2)| 
\leq c \Gamma_2(f_1) \Gamma_2(f_2)  
\]
for all $s_0 \in C$. 
Next, for  
$g\in U_1= C (U'')^m$ 
we write 
$g= s_0 g_1 g_2 \cdots g_m$ 
with $s_0\in C$ and $g_j\in U''$. 
Applying (\ref{ediffid}) shows that 
\[
(\partial_g f_1,f_2) = (\partial_{s_0} f_1,f_2) 
  + \sum_{j=1}^m
(\partial_{g_j} f_1, L((s_0 g_1 \cdots g_{j-1})^{-1}) f_2). 
\]
The lemma follows by applying the above estimate for $\partial_{s_0}$, 
the estimate of Lemma~\ref{luest}, and (\ref{enabla2}).  
\hfill $\Box$ 

\bigskip

In order to extend the estimate of Lemma~\ref{luoneest}, 
observe that $G_1 = \bigcup_{n=1}^{\infty} U_1^n$ 
since $U_1$ is a neighborhood of $e$ for the connected group $G_1$. 
We may therefore define a distance 
$\rho_1\colon G_1\to \Ni$ 
by 
\[
\rho_1(g) = \inf\{n\in\Ni\colon g\in U_1^n\}, 
\;\;\; 
g\in G_1. 
\]
From  
$G_1 \subseteq G$  
it follows easily that  
$\rho(g) \leq c \rho_1(g)$ for some constant $c>0$ 
and all 
$g\in G_1$.   
On the other hand, 
a theorem of Varopoulos \cite{Var20} asserts 
that any closed connected subgroup of a connected Lie group  
has at most exponential distance distortion.  
Applied to 
$G_1\subseteq G$, 
this means that there is 
$c'>0$ such that 
\begin{equation} 
\rho_1(g) \leq c' e^{c' \rho(g)},  \;\;\; g\in G_1.  
\label{edistort} 
\end{equation} 
We can now prove a version of (\ref{eg1gen}).

\begin{lemma} \label{lg1est} 
There exist $c,c'>0$ such that 
\[
|(\partial_g f_1,f_2)| 
\leq c e^{c' \rho(g)} \Gamma_2(f_1) \Gamma_2(f_2)  
\]
for all 
$g\in G_1$ and 
$f_1,f_2\in L^2$.  
\end{lemma} 
\proof\
Let 
$g\in G_1$ and set 
$n=\rho_1(g)$. 
Then 
$g=g_1 g_2 \cdots g_n$
for some  
$g_j\in U_1$.  
From (\ref{ediffid}) 
we have 
\[
(\partial_g f_1,f_2) 
= (\partial_{g_1}f_1,f_2) 
+ \sum_{j=2}^n (\partial_{g_j}f_1, L((g_1\cdots g_{j-1})^{-1})f_2). 
\]
Hence by Lemma~\ref{luoneest} and (\ref{enabla2}),  
\begin{eqnarray*} 
|(\partial_g f_1,f_2)| 
& \leq & c \Gamma_2(f_1) \Gamma_2(f_2)  
       + c \sum_{j=2}^n \Gamma_2(f_1)   
          \Gamma_2 (L((g_1\cdots g_{j-1})^{-1})f_2)   \\ 
& \leq & c \Gamma_2(f_1) \Gamma_2(f_2)  
        + c' \sum_{j=2}^n \rho(g_1 \ldots g_{j-1})
     \Gamma_2(f_1) \Gamma_2(f_2).    
\end{eqnarray*}   
Since 
\[
\rho(g_1 \cdots g_{j-1}) \leq \sum_{i=1}^{j-1} \rho(g_i) 
\leq c'' (j-1) \leq c'' n 
\]
where 
$c''= \sup \{\rho(g)\colon g\in U_1\}$, 
this yields an estimate 
\[
|(\partial_g f_1,f_2)| \leq c n^2 \Gamma_2(f_1) \Gamma_2(f_2). 
\] 
Recalling that 
$n=\rho_1(g)$, 
the lemma follows using  
(\ref{edistort}).  
\hfill $\Box$ 

\bigskip

\noindent {\bf Proof of Theorem~\ref{tgradupp1}}.  
Recall that 
$\pi_1\colon G\to G/G_1 \cong \Ri^q$;   
we identify 
$G/G_1$ with 
$\Ri^q$ 
so that 
$\pi_1(g) = (\pi_1^{(1)}(g), \ldots, \pi_1^{(q)}(g))\in\Ri^q$,   
$g\in G$.    

Let $\partial/(\partial y_1), \ldots, \partial/(\partial y_q)$ 
be the standard coordinate vector fields on $\Ri^q$, 
and fix 
$X_1, \ldots, X_q\in \gotg$ 
such that 
$(\pi_1)_* X_j = \partial/(\partial y_j)$ 
for all   
$j$. 
Let 
$\exp\colon \gotg\to G$ 
be the exponential map, 
and define the one-parameter subgroups 
$x_j\colon \Ri\to G$ 
by  
$x_j(t): = \exp(t X_j)$ 
for  
$t\in\Ri$, $j\in \{1, \ldots, q\}$.  

Given any 
$g\in G$, 
put 
$t_j: = \pi_1^{(j)}(g) \in\Ri$,  
$j\in \{1, \ldots, q\}$,  
and write   
\begin{equation} 
g = x_1(t_1) \cdots x_q(t_q) g' 
\label{egwrite} 
\end{equation} 
for some 
$g'\in G$.  
Because 
$\pi_1(g) = (t_1, \ldots, t_q) 
= \pi_1( x_1(t_1) \cdots x_q(t_q))$ 
and 
$G_1$ is the kernel of the homomorphism 
$\pi_1$,  
then 
$g' \in G_1$.  
In what follows, $c,c'$ will denote constants independent 
of $g\in G$. 
Since $\pi_1^{(j)}\colon G \to \Ri$ 
and  
$x_j\colon \Ri\to G$ 
are homomorphisms,  
it is straightforward to establish inequalities of form  
\begin{eqnarray} 
|t_j| & = & |\pi_1^{(j)}(g)| \leq c \rho(g),  \nonumber  \\   
\sup_{|s|\leq |t_j|} 
\rho(x_j(s)) & \leq & c \rho(g),
   \;\;\;  
\sup_{|u|\leq 1} \rho(x_j(u)) \leq c,  \nonumber \\
\rho(g') & \leq & c \rho(g),  
\label{emodineq1} 
\end{eqnarray} 
for all 
$g\in G$.  

Applying identity (\ref{ediffid2}) to
(\ref{egwrite})
shows that 
\begin{equation} 
\partial_g = \partial_{x_1(t_1)} + \cdots + \partial_{x_q(t_q)} 
    + \partial_{g'} + W_g,  
\label{epartgdec} 
\end{equation} 
where the operator $W_g$ is 
a finite sum of terms each of the form 
$\partial_{h_1} \cdots \partial_{h_k}$ 
with $k\geq 2$ and 
$h_1, \ldots, h_k  \in \{x_1(t_1), \ldots, x_q(t_q), g'\}$.   
By writing 
\begin{eqnarray*} 
|(\partial_{h_1} \ldots \partial_{h_k} f_1,f_2)| 
 & = & 
  |(\partial_{h_2} \ldots \partial_{h_k} f_1, \partial_{h_1^{-1}}f_2)| 
  \\
 & \leq & \| \partial_{h_2} \ldots \partial_{h_k} f_1 \|_2 
          \| \partial_{h_1^{-1}}f_2\|_2 
 \leq 2^{k-2} \| \partial_{h_k} f_1\|_2 \| \partial_{h_1^{-1}}f_2\|_2, 
\end{eqnarray*} 
one sees 
from (\ref{enabla1}) and (\ref{emodineq1}) 
that 
\begin{equation} 
|(W_g f_1,f_2)| \leq c \rho(g)^2 \Gamma_2(f_1) \Gamma_2(f_2)  
\label{eauxe1}  
\end{equation} 
for all 
$f_1,f_2\in L^2$.  
Also, 
since 
$g'\in G_1$,  
Lemma~\ref{lg1est} and (\ref{emodineq1}) 
provide an estimate 
\begin{equation} 
|(\partial_{g'} f_1,f_2)| \leq c e^{c\rho(g')} 
\Gamma_2(f_1) \Gamma_2(f_2) 
\leq c e^{c' \rho(g)} \Gamma_2(f_1) \Gamma_2(f_2) 
\label{eauxe2} 
\end{equation} 
for 
$f_1,f_2\in L^2$. 
Next, 
the following lemma allows us to compare 
$\partial_{x_j(t_j)}$ and 
$t_j \partial_{x_j(1)}$.

\begin{lemma} \label{lonegp} 
Let 
$(R_t)_{t\in\Ri} 
\subseteq {\cal L}(X)$ 
be a strongly continuous one-parameter 
group of bounded linear operators in a Banach space 
$X$. 
Set 
$D_t:= R_t-I$ 
for 
$t\in\Ri$.  
Then 
\[
D_t - t D_1  = -\int_0^1 ds\, D_s D_t  + \int_0^t ds\, D_s D_1    
\]
where the integrals converge with respect to the strong operator 
topology.    
\end{lemma} 
\noindent {\bf Proof of Lemma~\ref{lonegp}.}  
For 
$t\in\Ri$ 
set 
\begin{eqnarray*} 
P_t  & := & \int_0^t ds\, (R_{s+1}-R_s)   \\
  & = & \int_1^{t+1} ds\, R_s  - \int_0^t ds\, R_s    
    = \int_0^1 ds\,  (R_{s+t}-R_s), 
\end{eqnarray*} 
and observe, using the group property 
$R_u R_v = R_{u+v}$,  
that  
\[
D_t-P_t = \int_0^1 ds\, (D_t -(R_{s+t}-R_s)) 
        = -\int_0^1 ds\, D_s D_t,  
\]
\[
P_t - t D_1 = \int_0^t ds\, ((R_{s+1}-R_s)-D_1) 
           = \int_0^t ds\, D_s D_1.  
\]
The lemma follows.    
\hfill $\Box$ 

\bigskip 

Setting 
$R_t = L(x_j(t))$ 
and  
$D_t = \partial_{x_j(t)}$ 
in Lemma~\ref{lonegp}
gives 
\[
\partial_{x_j(t_j)} - t_j \partial_{x_j(1)} 
 = -\int_0^1 ds\, \partial_{x_j(s)} \partial_{x_j(t_j)} 
    + \int_0^{t_j} ds\, \partial_{x_j(s)} \partial_{x_j(1)}. 
\]
From this identity and  
(\ref{enabla1}),   
(\ref{emodineq1}),
one easily sees that 
\begin{equation}
| ((\partial_{x_j(t_j)} - t_j \partial_{x_j(1)})f_1,f_2)| 
\leq c \rho(g)^2  \Gamma_2(f_1) \Gamma_2(f_2)
\label{eauxe3}
\end{equation} 
for $f_1,f_2\in L^2$. 

By combining (\ref{eauxe1}), (\ref{eauxe2}), (\ref{eauxe3})  
with (\ref{epartgdec})    
and recalling that 
$t_j = \pi_1^{(j)}(g)$,  
we conclude that 
\begin{equation} 
\partial_g = \sum_{j=1}^q \pi_1^{(j)}(g) \partial_{x_j(1)} 
  + V_g,  
\label{eveedef} 
\end{equation} 
where the linear operator 
$V_g\in {\cal L}(L^2)$ satisfies  
\begin{equation} 
| (V_g f_1,f_2)| \leq 
c e^{c\rho(g)} \Gamma_2(f_1) \Gamma_2(f_2) 
\label{evesttwo} 
\end{equation} 
for all 
$g\in G$, 
$f_1,f_2\in L^2$.    
But since $K$ is centered,  
it follows from (\ref{elapdef}) and 
(\ref{edefcent})
that  
\[
(I-T)f = - \int_G dg\, K(g) \partial_g f = -\int_G dg\, K(g) V_g f 
\]
for all 
$f\in L^2$. 
Therefore, by (\ref{evesttwo}) and 
Assumption~(ii),  
\begin{eqnarray*} 
|((I-T)f,f)| 
& \leq & \int_G dg\, K(g) |(V_g f,f)|     \\ 
& \leq & \int_G dg\, K(g) ce^{c \rho(g)} \Gamma_2(f)^2 
\leq c' \Gamma_2(f)^2  
\end{eqnarray*}  
for all 
$f\in L^2$.  
This ends the proofs of Theorems~\ref{tgradupp1} 
and \ref{tcentanal1}.     
\hfill $\Box$ 

\bigskip 

\noindent {\bf Remarks.} 

(a) Assumption~(ii) of Theorem~\ref{tcentanal1} 
can be relaxed in case $G$ has polynomial 
volume growth.   
Indeed, in this case,  
by a result of 
\cite{Var20} (see also \cite{Boun}) 
any  closed subgroup of $G$ has polynomial 
distance distortion, 
which means that there exists 
$N\geq 1$ with 
\begin{equation} 
\rho_1(g) \leq c \rho(g)^N 
\label{epolydist} 
\end{equation} 
for all 
$g\in G_1$.  
Then, 
one can prove that 
$T$ is analytic
for any centered density $K$
satisfying Assumption~(i) and the condition 
\[ 
\int_G dg\, K(g) \rho(g)^{2N} < \infty. 
\]
This result follows by an easy adaption of the above proofs;   
one uses (\ref{epolydist}) to obtain   
Lemma~\ref{lg1est} in the form 
\[
|(\partial_g f_1,f_2)| 
\leq c\rho(g)^{2N} \Gamma_2(f_1) 
\Gamma_2(f_2)
\] 
for all 
$f_i\in L^2$ and  
$g\in G_1$.

\bigskip

(b) Let us describe an 
$L^p$ version of the $L^2$ estimate of 
Theorem~\ref{tgradupp1}.   
  
Fix 
$p\in [1,\infty]$ 
with 
$p^{-1} + (p')^{-1}=1$      
and consider the pairing 
$(f_1,f_2) : = \int_G dg\, f_1 \overline{f_2}$ 
for 
$f_1\in L^p$, 
$f_2\in L^{p'}$.  
By modifying the arguments leading to 
(\ref{evesttwo}), 
one obtains 
\begin{equation}  
|(V_g f_1, f_2)| \leq c e^{c\rho(g)} \Gamma_p(f_1) \Gamma_{p'}(f_2) 
\label{evestp} 
\end{equation} 
for all 
$g\in G$, 
$f_1\in L^p$,  
$f_2\in L^{p'}$, 
where 
$V_g$ is as in (\ref{eveedef}). 
To see this,    
one uses H\"{o}lder's inequality 
in place of the Cauchy-Schwarz inequality   
in the proofs of Lemmas~\ref{luest}, \ref{luoneest}, 
\ref{lg1est} and Theorem~\ref{tgradupp1}.   

When $K$ is centered and satisfies Assumption~(ii), 
estimate~(\ref{evestp}) yields that 
\begin{equation} 
|((I-T)f_1,f_2)| 
= \left| \int_G dg\, K(g) (V(g)f_1,f_2) 
\right| 
\leq c \Gamma_p(f_1) \Gamma_{p'}(f_2) 
\label{euppgradp} 
\end{equation} 
for all $f_1\in L^p$ and $f_2\in L^{p'}$. 
In particular, 
with 
$p=\infty$ we get that    
\begin{equation} 
|((I-T)\psi, f)| \leq c \Gamma_1(f) 
\label{euppgradpsi} 
\end{equation} 
for all 
$\psi\in {\cal E} \cap L^{\infty}$ and 
$f\in L^1$.  
Because the right side of (\ref{euppgradpsi}) 
is independent of 
$\|\psi\|_{\infty}$,  
an easy approximation argument then gives 
(\ref{euppgradpsi}) 
for all 
$\psi\in {\cal E}$.

The estimate~(\ref{euppgradpsi}) for centered $K$ 
will be crucial for the proof of  Theorem~\ref{tcentgauss1}.

\section{Proof of Theorem~\ref{tcentgauss1}}  \label{s4}  

For the proof of Theorem~\ref{tcentgauss1}, 
in this section 
$K$ will denote a sub-Gaussian density 
on the connected Lie group 
$G$. 
Let us fix 
$\psi\in {\cal E}$ 
and set 
$T_{\lambda}: = e^{\lambda\psi} T e^{-\lambda\psi}$ 
for 
$\lambda\in\Ri$.   
The constants in the estimates obtained below do not, however, depend on 
the particular choice of  
$\psi\in {\cal E}$.   

Note that, by 
(\ref{enabla1}), 
we have 
$\|\partial_g \psi\|_{\infty} \leq \rho(g)$
for all 
$g\in G$.

We begin with a simple estimate of 
$T_{\lambda}-T$  
which is valid whether or not $K$ is centered. 

\begin{lemma} \label{lcrudest} 
There exist constants $c,\omega>0$ such that 
$\|T_{\lambda}-T\|_{p\to p} 
\leq c|\lambda| e^{\omega\lambda^2}$ 
for all 
$\lambda\in\Ri$ 
and 
$p\in [1,\infty]$.     
\end{lemma} 
\proof\
Note the identity 
\[
e^{\lambda\psi} L(g) (e^{-\lambda\psi}f) 
 = e^{-\lambda\partial_g \psi} L(g)f 
\]
for a function $f\colon G\to \Ci$ 
and $g\in G$.  
Since 
$T=\int_G dg\, K(g) L(g)$, 
then 
\begin{equation} 
(T_{\lambda}-T)f = \int_G dg\, K(g) 
[e^{-\lambda\partial_g \psi}-1] L(g)f. 
\label{ediffexp} 
\end{equation} 
The inequalities  
$|e^s-1| \leq |s| e^{|s|}$, 
$s\in\Ri$, 
and 
$\|\partial_g \psi\|_{\infty}\leq \rho(g)$ 
imply that   
\[
\| e^{-\lambda\partial_g \psi} -1 \|_{\infty} 
\leq |\lambda|\rho(g) e^{|\lambda|\rho(g)} 
\]
for all $g\in G$, $\lambda\in\Ri$. 
Note the elementary estimate 
$|\lambda|\rho(g) \leq \varepsilon\rho(g)^2 + \varepsilon^{-1}\lambda^2$ 
for all 
$\varepsilon>0$.   
Fixing 
$\varepsilon<b$,  
where 
$b$ is as in (\ref{ekgauss1}),  
we get    
\begin{eqnarray*} 
\|(T_{\lambda}-T)f\|_p 
& \leq & 
|\lambda| \int_G dg\, K(g) \rho(g) e^{|\lambda| \rho(g)} \|f\|_p  
  \\\
& \leq & c|\lambda| e^{\varepsilon^{-1} \lambda^2}  
\int_G dg\, e^{-(b-\varepsilon)\rho(g)^2} \rho(g)  
\|f\|_p 
\leq c' |\lambda| e^{\varepsilon^{-1}\lambda^2} \|f\|_p. 
\end{eqnarray*} 
The lemma follows. 
\hfill $\Box$ 

\bigskip 

Since 
$\|T\|_{2\to 2} \leq 1$,  
Lemma~\ref{lcrudest} implies that 
$\|T_{\lambda}\|_{2\to 2} = 1 + O(|\lambda|)$ 
for small $|\lambda|$. 
To prove Theorem~\ref{tcentgauss1}, 
our main task will be to obtain the improved estimate 
\begin{equation} 
\|T_{\lambda}\|_{2\to 2} \leq 1 + c\lambda^2 
\label{elambdasq} 
\end{equation} 
for all $|\lambda|\leq 1$, 
when $K$ is centered.   
Indeed, given (\ref{elambdasq}) 
and Lemma~\ref{lcrudest} it follows for some $\omega'>0$ 
that 
$\|T_{\lambda}\|_{2\to 2} \leq e^{\omega' \lambda^2}$ 
for all $\lambda\in\Ri$.  
This implies 
$\|T_{\lambda}^n\|_{2\to 2} \leq e^{\omega' \lambda^2 n}$  
for all $n\in\Ni$, 
which is the estimate of Theorem~\ref{tcentgauss1}. 

Define quadratic forms $Q$ and $Q_{\lambda}$, 
$\lambda\in\Ri$, 
by 
\[
Q(f) = \|f\|_2^2 - \|Tf\|_2^2, 
\;\;\; 
Q_{\lambda}(f) = \|f\|_2^2 - \|T_{\lambda}f\|_2^2, 
\]
for all $f\in L^2$, so that 
$Q_0=Q$.  
To obtain (\ref{elambdasq}) we require the 
perturbation estimate of the next proposition.  
Note that (\ref{euppgradpsi}) is essential 
for the proof of the proposition.   

\begin{prop} \label{ppert}  
\noindent {\rm (I)}$\;$  
There exists 
$c>1$ 
such that 
\[
c^{-1} \Gamma_2(f)^2 \leq Q(f) \leq c \Gamma_2(f)^2 
\]
for all $f\in L^2$. 

\noindent {\rm (II)}$\;$  
If $K$ is centered, then 
there exists $c'>0$ such that 
\[
|Q(f)-Q_{\lambda}(f)| \leq \varepsilon Q(f) + c'(1+\varepsilon^{-1}) 
 \lambda^2 \|f\|_2^2 
\]
for all 
$\varepsilon>0$, 
$f\in L^2$ and 
$|\lambda|\leq 1$.  
\end{prop} 
\proof\
First 
observe that 
$Q(f) = (f,(I-T^* T)f)$ 
and 
$T^* Tf = \widetilde{K} * f$, $f\in L^2$, 
where 
$\widetilde{K}$ is the density on $G$ defined 
by 
$\widetilde{K}: = K^* * K$ 
with 
$K^*(g): = \Delta(g^{-1}) K(g^{-1})$. 
Then part~(I) follows easily by applying Lemma~\ref{lrealest} 
to $T^* T$ in place of $T$. 

To prove part~(II), it suffices to prove an  estimate of form  
\begin{equation} 
| Q(f)-Q_{\lambda}(f)| \leq c|\lambda| \Gamma_2(f) \|f\|_2 
 + c \lambda^2 \|f\|_2^2 
\label{eqpert2} 
\end{equation} 
for all $f\in L^2$, $|\lambda|\leq 1$. 
For then part~(II) follows from part~(I) and the 
elementary estimate 
$|\lambda| Q(f)^{1/2} \|f\|_2 \leq \varepsilon Q(f)   
+ \varepsilon^{-1}  \lambda^2 \|f\|_2^2$ 
for all 
$\varepsilon>0$. 
 
We begin the proof of (\ref{eqpert2}) 
by writing 
\begin{eqnarray*} 
Q(f) -Q_{\lambda}(f) & = & \|T_{\lambda}f\|_2^2 - \|Tf\|_2^2  \\
& = & (Tf, (T_{\lambda}-T)f) + ((T_{\lambda}-T)f,Tf)   \\
    & & {}+ ((T_{\lambda}-T)f, (T_{\lambda}-T)f)    \\ 
    & = & 2 \RRe ((T_{\lambda}-T)f,Tf) + \|(T_{\lambda}-T)f\|_2^2 \\
    & = & 2 \RRe ((T_{\lambda}-T)f,f) 
       + \|(T_{\lambda}-T)f\|_2^2 \\ 
    & & {} - 2\RRe ((T_{\lambda}-T)f,(I-T)f).    
\end{eqnarray*}   
By Lemma~\ref{lcrudest},  
the term 
$\|(T_{\lambda}-T)f\|_2^2$ 
on the right side is estimated 
by 
$c \lambda^2 \|f\|_2^2$ for 
$|\lambda|\leq 1$. 
Also, from (\ref{elapdef}) we have 
$\|(I-T)f\|_2 \leq c \Gamma_2(f)$ 
so that 
\[
|((T_{\lambda}-T)f,(I-T)f)| \leq c |\lambda| \|f\|_2 \Gamma_2(f)
\]
for all $f\in L^2$.  
To handle the remaining term 
$((T_{\lambda}-T)f,f)$, 
by (\ref{ediffexp})  
we may write 
\begin{eqnarray*} 
((T_{\lambda}-T)f,f) 
& = & \int_G dg\, K(g) ([e^{-\lambda\partial_g \psi}-1] L(g)f,f)  
     \\
& = & \int_G dg\, K(g) 
([e^{-\lambda\partial_g \psi}-1+\lambda \partial_g\psi]
L(g)f,f)   \\ 
  & & {}  - \lambda  \int_G dg\, K(g) ((\partial_g \psi) \partial_g f, f)  
            \\
  & & {} - \lambda  \int_G dg\, K(g) ((\partial_g \psi) f,f)  
            \\
  & = & I_1 + I_2 + I_3, 
\end{eqnarray*}    
using that 
$L(g)f = \partial_g f + f$.   
We will bound each of the terms 
$I_1$, 
$I_2$, 
$I_3$.   
First, 
since 
$|e^s- 1 -s| \leq s^2 e^{|s|}$ for $s\in\Ri$ 
and
$\|\partial_g \psi\|_{\infty}\leq \rho(g)$,
$\|L(g)f\|_2=\|f\|_2$, 
one obtains 
\[
|([e^{-\lambda\partial_g \psi}-1 + \lambda\partial_g \psi] L(g)f, f)| 
\leq  
\lambda^2 \rho(g)^2 e^{|\lambda|\rho(g)} \|f\|_2^2 
\]
for all $g\in G$.    
It easily follows that 
$|I_1| \leq c \lambda^2 \|f\|_2^2$ 
for all 
$|\lambda|\leq 1$.  

Next, observing that 
\[
|((\partial_g\psi)\partial_g f, f)| 
\leq \|\partial_g \psi\|_{\infty} 
\|\partial_g f\|_2 \|f\|_2 
\leq \rho(g)^2 \Gamma_2(f) \|f\|_2, 
\]
we deduce that 
$|I_2| \leq c |\lambda| \Gamma_2(f) \|f\|_2$ 
for all 
$|\lambda|\leq 1$.  

Finally, from (\ref{elapdef}) we can write 
\[
I_3 = \lambda (((I-T)\psi)f, f) 
    = \lambda ((I-T)\psi,  f \overline{f}).   
\]
Since $K$ is centered 
(this is the only place in the argument 
where centeredness is used), 
(\ref{euppgradpsi}) gives 
\[
|((I-T)\psi, f \overline{f})| 
\leq c \Gamma_1 (f \overline{f}) 
\leq 2c \Gamma_2(f) \|f\|_2  
\]
for all $f\in L^2$, 
where the second inequality follows easily from the identity  
$\partial_u(f_1 f_2) 
= (\partial_u f_1) f_2 + (L(u)f_1)(\partial_u f_2)$
for 
$u\in U$, $f_1,f_2\in L^2$.  
Thus 
$|I_3| \leq c'|\lambda| \Gamma_2(f) \|f\|_2$ 
for all 
$|\lambda|\leq 1$. 
Collecting the above estimates  
yields (\ref{eqpert2}),  
and Proposition~\ref{ppert} is proved.  
\hfill $\Box$ 

\bigskip

\noindent {\bf Remark.} 
The proof of Proposition~\ref{ppert} also yields the following 
perturbation estimate of interest:  
for centered 
$K$, 
one has 
\begin{equation} 
|((T_{\lambda}-T)f,f)| 
\leq \varepsilon \Gamma_2(f)^2 
+ c(1+\varepsilon^{-1})\lambda^2 \|f\|_2^2 
\label{epert3} 
\end{equation} 
for all 
$\varepsilon>0$, $|\lambda|\leq 1$ 
and 
$f\in L^2$. 
Actually, 
our argument above shows that 
(\ref{epert3}) is essentially equivalent with 
the estimate of Proposition~\ref{ppert}, Part~(II).  

\bigskip

To continue with the proof of (\ref{elambdasq}),  
choose 
$\varepsilon=1/2$ in the estimate of Proposition~\ref{ppert}, part~(II). 
We find,  
for some 
$c>0$,  
that  
\begin{eqnarray*} 
\|f\|_2^2 - \|T_{\lambda}f\|_2^2 
& = & Q_{\lambda}(f) \geq Q(f) - |Q(f)-Q_{\lambda}(f)| \\ 
& \geq & 2^{-1} Q(f) - c\lambda^2 \|f\|_2^2 
\geq -c\lambda^2 \|f\|_2^2 
\end{eqnarray*} 
for all $f\in L^2$ and $|\lambda|\leq 1$. 
Thus 
$\|T_{\lambda}f\|_2^2 \leq (1+c\lambda^2)\|f\|_2^2$ 
and 
$\|T_{\lambda}\|_{2\to 2} \leq (1+c\lambda^2)^{1/2} 
\leq 1+c\lambda^2$ 
for 
$|\lambda|\leq 1$.  
This proves (\ref{elambdasq}) 
and completes the proof of 
Theorem~\ref{tcentgauss1}.   
\hfill $\Box$

\section{Proof of Theorem~\ref{tcentgset}} \label{s5}

Theorem~\ref{tcentgset} essentially follows from 
Theorem~\ref{tcentgauss1}
and the equivalence between two different forms 
of $L^p$ off-diagonal estimates.   

On any locally compact, compactly generated group $G$, 
we can give the following precise statement of this equivalence. 
(If 
$E$ or $F$ is the empty set,  
adopt the convention that 
$d(E,F)= +\infty$ and 
$e^{-b d(E,F)^2}= 0$ for 
$b>0$.)     

\begin{prop} \label{pgexp}
Fix $p\in [1,\infty]$, 
let $S\in {\cal L}(L^p)$ and let $n\in\Ni$.  

\noindent {\rm (I)}$\;$ 
If there exist constants 
$a,\omega>0$ such that 
\[
\| e^{\lambda\psi} S e^{-\lambda\psi} \|_{p\to p} 
 \leq a e^{\omega \lambda^2 n} 
\]
for all $\lambda\in\Ri$ and $\psi\in {\cal E}$, 
then there exist $c,b>0$ depending on $\omega$ but 
not on $n$ or $a$, 
such that 
\[
\| \chi_E S \chi_F \|_{p\to p} \leq ca e^{-b d(E,F)^2/n} 
\]
for all Borel measurable subsets $E,F\subseteq G$.  

\noindent {\rm (II)}$\;$ 
Conversely, 
if there exist 
$a,b>0$ such that 
\[
\| \chi_E S \chi_F \|_{p\to p} \leq a e^{-bd(E,F)^2/n} 
\]
for all Borel measurable $E,F\subseteq G$,  
then there exist $c,\omega>0$ depending on $b$
but not on $n$ or $a$, 
such that 
\[
\| e^{\lambda \psi} S e^{-\lambda \psi} \|_{p\to p} 
\leq ca e^{\omega \lambda^2 n} 
\]
for all $\lambda\in\Ri$ and $\psi\in {\cal E}$. 
\end{prop} 

The case 
$p=2$ of Theorem~\ref{tcentgset} is 
an immediate consequence of Theorem~\ref{tcentgauss1} 
and Proposition~\ref{pgexp}. 
The general case 
$p\in (1,\infty)$ of Theorem~\ref{tcentgset} 
then follows by interpolation between the case 
$p=2$ 
and the obvious estimates 
$\| \chi_E T^n \chi_F \|_{1\to 1} \leq \|T^n \|_{1\to 1} \leq 1$, 
$\| \chi_E T^n \chi_F \|_{\infty\to \infty} \leq 1$.  

Thus it only remains to prove 
Proposition~\ref{pgexp}, 
and the rest of this section is occupied with this proof.   

The proof of part~(I) 
is a variation of standard ideas  
(see for example \cite{Dav10}.)    
To prove part~(I), 
let $E,F$ be non-empty 
with 
$d(E,F)\geq 2$ 
(if 
$d(E,F)=1$ the desired estimate 
is trivial since  
$\|\chi_E S \chi_F \|_{p\to p} \leq \|S\|_{p\to p} \leq a$.)
We may assume that 
$a=1$ 
(otherwise, 
replace 
$S$ by 
$a^{-1}S$). 
Define 
$\psi_F\colon G\to \Ni$ by 
$\psi_F(g) = d(\{g\},F) 
= \inf\{\rho(gh^{-1})\colon h\in F\}$.  
The bound 
\[
|\psi_F(h) - \psi_F(gh)| \leq \rho(g), 
\;\;\; 
g,h\in G,  
\]
implies that 
$\psi_F \in {\cal E}$.  
For 
$g\in E$ we have $\psi_F(g) \geq d(E,F)$, 
while for 
$g\in F$, 
$\psi_F(g)=1\leq 2^{-1}d(E,F)$. 
For  
$\lambda\geq 0$, 
then  
\[
\| e^{-\lambda \psi_F} \chi_E \|_{\infty} \leq e^{-\lambda d(E,F)}, 
\;\;\; 
\| e^{\lambda\psi_F} \chi_F \|_{\infty} \leq e^{2^{-1}\lambda d(E,F)}. 
\]
Applying the hypothesis of (I) gives 
\begin{eqnarray*} 
\| \chi_E S \chi_F \|_{p\to p} 
& \leq & \| e^{-\lambda \psi_F} \chi_E \|_{\infty} 
         \| e^{\lambda \psi_F} S e^{-\lambda \psi_F} 
        \|_{p\to p} 
         \| e^{\lambda \psi_F} \chi_F \|_{\infty}    \\
& \leq & e^{\omega \lambda^2 n - 2^{-1} \lambda d(E,F)}
\end{eqnarray*}
for all 
$\lambda\geq 0$.   
Choosing $\lambda$ to be a small constant multiple 
of 
$d(E,F)/n$ proves part~(I).    

For the proof of part~(II) we need  the following lemma, 
which actually  
holds on any measure space.  

\begin{lemma} \label{lorthog} 
Let $p,p'\in [1,\infty]$ 
with 
$p^{-1} + (p')^{-1} =1$.  
Let 
$\{E_k\}$ be a finite or countable sequence of  
measurable subsets of $G$  
such that 
$G = \bigcup_k E_k$.  
Then any linear operator $S$ in $L^p$ satisfies 
\[
\| S\|_{p\to p}
\leq \left( \sup_l \sum_k \|\chi_{E_k} S \chi_{E_l}\|_{p\to p} 
       \right)^{1/p} 
     \left( \sup_k \sum_l \|\chi_{E_k} S \chi_{E_l}\|_{p\to p} 
       \right)^{1/p'}. 
\]
\end{lemma} 
\proof\
We may assume that the $E_k$ are pairwise disjoint 
(otherwise, we may replace $\{E_k\}$ with a pairwise disjoint sequence 
$\{E_k'\}$ such that 
$E_k'\subseteq E_k$ 
and 
$\bigcup_k E_k' = \bigcup_k E_k = G$.)     
Since the cases 
$p=1$ or $p=\infty$ are straightforward, 
we will also assume that 
$p\in (1,\infty)$.  
Suppose 
$\|f_1\|_p=1$ and  
$\|f_2\|_{p'}=1$.  
Set 
$f_j^{(k)}: = \chi_{E_k} f_j$  
for 
$j=1,2$,  
and observe that   
\[
(S f_1, f_2) = \sum_{k,l} (\chi_{E_k} S \chi_{E_l} f_1^{(l)}, 
                                         f_2^{(k)}). 
\]
By applications of H\"{o}lder's inequality,  
we obtain  
\begin{eqnarray*} 
|(Sf_1,f_2)| & \leq & \sum_{k,l} \|\chi_{E_k} S \chi_{E_l} \|_{p\to p} 
                       \| f_1^{(l)}\|_p \|f_2^{(k)}\|_{p'}     \\ 
& = & \sum_{k,l} (\|\chi_{E_k} S \chi_{E_l}\|_{p\to p}^{1/p} 
              \|f_1^{(l)}\|_p)                      
                 (\|\chi_{E_k} S \chi_{E_l}\|_{p\to p}^{1/p'}
              \|f_2^{(k)}\|_{p'})             \\
& \leq & \left( \sum_{k,l} \|\chi_{E_k} S \chi_{E_l} \|_{p\to p} 
                 \|f_1^{(l)}\|_p^p \right)^{1/p}   \\
& &  {} \times        \left( \sum_{k,l} \|\chi_{E_k} S \chi_{E_l} \|_{p\to p} 
                 \|f_2^{(k)}\|_{p'}^{p'} \right)^{1/p'}   \\
& \leq & \left( \sup_l \sum_k \|\chi_{E_k} S \chi_{E_l}\|_{p\to p} 
       \right)^{1/p} 
     \left( \sup_k \sum_l \|\chi_{E_k} S \chi_{E_l}\|_{p\to p} 
       \right)^{1/p'},  
\end{eqnarray*} 
where the last step used that   
$\sum_l \|f_1^{(l)}\|_p^p = \|f_1\|_p^p=1$ 
and 
$\sum_k \|f_2^{(k)}\|_{p'}^{p'} 
= \|f_2\|_{p'}^{p'} =1$. 
The lemma follows. 
\hfill $\Box$ 

\bigskip 

To prove part~(II), 
we may assume that 
$a=1$. 
Fix 
$\psi\in {\cal E}$, 
define 
\[
E_k : = \{g\in G\colon kn^{1/2} \leq \psi(g) < (k+1)n^{1/2}\}, 
\;\;\; 
k\in\Zi,     
\]
and observe that $G$ is the disjoint union of the 
$E_k$, $k\in\Zi$.  
The conclusion of part~(II) will 
follow from  Lemma~\ref{lorthog},  
if we can prove a bound  
\begin{equation} 
\| \chi_{E_k} e^{\lambda\psi} S e^{-\lambda\psi} \chi_{E_l}\|_{p\to p} 
\leq c e^{-b'|k-l|^2} e^{\omega\lambda^2 n} 
\label{epertor1} 
\end{equation} 
for all 
$\lambda\in\Ri$, $k,l\in \Zi$, 
where 
$c,b',\omega$ are positive constants independent of 
$n$ and $\psi$.
We will prove (\ref{epertor1}) for 
$\lambda\geq 0$ 
(the case 
$\lambda\leq 0$ being similar).  
Since 
\[
\|\chi_{E_k} e^{\lambda \psi}\|_{\infty} 
\leq e^{\lambda(k+1)n^{1/2}}, 
\;\;\; 
\|\chi_{E_l} e^{-\lambda\psi}\|_{\infty} 
\leq e^{-\lambda l n^{1/2}},  
\]
we get 
\begin{eqnarray} 
\| \chi_{E_k} e^{\lambda\psi} S e^{-\lambda\psi} \chi_{E_l}\|_{p\to p} 
& \leq & e^{\lambda (|k-l|+1) n^{1/2}} 
       \|\chi_{E_k} S \chi_{E_l}\|_{p\to p}    \nonumber \\
& \leq & e^{\lambda (|k-l|+1) n^{1/2}} e^{-b d(E_k,E_l)^2/n}
\label{epertor2} 
\end{eqnarray}  
for all 
$\lambda\geq 0$ and $k,l\in\Zi$. 
In case $|k-l|\leq 1$,  it follows that 
\[
\|\chi_{E_k} e^{\lambda\psi} S e^{-\lambda\psi} \chi_{E_l}\|_{p\to p} 
\leq e^{2\lambda n^{1/2}} 
\leq e^{1+\lambda^2 n},  
\]
which implies (\ref{epertor1}) in this case. 

In the remaining case where 
$|k-l|\geq 2$,  
since 
$\psi\in {\cal E}$ 
we have 
\[
\rho(gh^{-1}) \geq |\psi(g) - \psi(h)| 
\geq (|k-l|-1)n^{1/2} 
\]
whenever 
$g\in E_k$ and $h\in E_l$. 
Hence 
\[
d(E_k,E_l) \geq (|k-l|-1)n^{1/2} \geq 2^{-1} |k-l|n^{1/2} 
\]
and,  by  
(\ref{epertor2}),  
\[
\|\chi_{E_k} e^{\lambda\psi} S e^{-\lambda\psi} \chi_{E_l}\|_{p\to p} 
\leq e^{2\lambda |k-l|n^{1/2}} e^{-4^{-1} b |k-l|^2} 
\]
for all  
$|k-l|\geq 2$ and 
$\lambda\geq 0$.  
By writing 
$2\lambda |k-l|n^{1/2} \leq \varepsilon |k-l|^2 
 + \varepsilon^{-1} \lambda^2 n$
for a small constant 
$\varepsilon>0$, 
we deduce a bound (\ref{epertor1}) 
for  
$|k-l|\geq 2$. 
The proof of Proposition~\ref{pgexp} is complete. 
\hfill $\Box$

\section{Proof of Theorem~\ref{tamen}}   \label{s6}

The proof of Theorem~\ref{tamen} 
is very similar to the proof of the analogous 
theorem for discrete groups given  in 
\cite[Theorem~1.9]{Dun20}. 
For this reason,  
we will only outline the arguments.  

The implications (I)$\Rightarrow$(III)
and (III)$\Rightarrow$(II)  
follow from the proofs of 
Theorems~\ref{tgradupp1} 
and \ref{tcentanal1}.  
The implications 
(I)$\Rightarrow$(V)$\Rightarrow$(IV)   
follow from the proof of Theorem~\ref{tcentgauss1}. 
We also have 
(I)$\Rightarrow$(VI)  
by Theorem~\ref{tcentgset}, 
and 
(VI)$\Rightarrow$(IV) by Proposition~\ref{pgexp}.  

It remains to show 
(II)$\Rightarrow$(I) and 
(IV)$\Rightarrow$(I).   
Suppose that 
$K$ is a non-centered, sub-Gaussian density on the amenable Lie group 
$G$.   
Fix a 
$j\in \{1, \ldots, q\}$ 
such that 
$\int_G dg\, K(g)\pi_1^{(j)}(g) \neq 0$.  
Let 
$P$
be the image of 
$K$ 
under the homomorphism 
$\pi_1^{(j)}\colon G\to \Ri$;  
that is, 
$P\colon \Ri\to [0,\infty)$ 
is the density on 
$\Ri$ such that 
\[ 
\int_\Ri dx\, P(x) F(x) = \int_G dg\, K(g) F(\pi_1^{(j)}(g)) 
\]
for all 
$F\in L^{\infty}(\Ri)$.
Then
$\int_\Ri dx\, x P(x)\neq 0$, 
in other words, 
$P$ is not centered.    
Let 
$\psi(x):=x$, 
$x\in\Ri$, 
and set 
\[
Sf: = P*f, 
\;\;\;   
S_{\lambda}f : = e^{\lambda\psi} S (e^{-\lambda\psi}f) 
= (e^{\lambda \psi}P) * f  
\]
for all 
$f\in L^p(\Ri)$ and 
$\lambda\in\Ri$.  
Using the Fourier theory of 
$L^2(\Ri)$,  
it is straightforward to show that  
$S$ is not analytic in 
$L^2(\Ri)$
and that the estimate 
$\|S_{\lambda}^n\|_{2\to 2} \leq c e^{\omega \lambda^2 n}$, 
$n\in\Ni$, 
$\lambda\in [-1,1]$, 
cannot hold.  
(Compare, for example, 
\cite[Lemma~6.2]{Dun20} for similar results for 
non-centered densities on 
$\Zi$.)   
But the transference theorem for convolution operators on amenable groups 
(see \cite[Theorem~2.4]{CW})   
implies that 
\[
\|(I-S)S^n\|_{2\to 2} \leq \|(I-T)T^n\|_{2\to 2}, 
\;\;\; 
\|S_{\lambda}^n \|_{2\to 2} \leq 
 \| e^{\lambda \pi_1^{(j)}} T^n e^{-\lambda\pi_1^{(j)}}\|_{2\to 2} 
\]
for all $n$ and $\lambda$.  
Consequently, 
Conditions~(II) and (IV) fail 
when $G$ is amenable and 
$K$ is not centered. 
\hfill $\Box$

\section{Further proofs} \label{s7} 

In this section we describe  
the proofs of Theorems~\ref{tcentgauss2} to \ref{tunifanal}.   
We will be brief, since some of the arguments occur elsewhere 
in a similar form.

Theorem~\ref{tcentgauss2} follows from 
Theorem~\ref{tcentgauss1} and (\ref{ekgauss1}) 
via known arguments: 
compare for example  
the proof of \cite[Theorem~1.2]{Dun15},  
or \cite[pp.126-127]{VSC}.  
The idea is to observe that, 
setting 
$\widetilde{T}_{\lambda}: 
= \Delta^{1/2} e^{\lambda\rho} T e^{-\lambda\rho} \Delta^{-1/2}$, 
the operators 
$\widetilde{T}_{\lambda}^n$ 
have integral kernels with respect to 
the {\it right} Haar measure 
$d\widehat{g} = \Delta(g^{-1}) dg$ 
given by 
\[
\widetilde{K}_{n,\lambda}(g,h) 
:= e^{\lambda \rho(g)} \Delta^{1/2}(g) K^{(n)}(gh^{-1}) 
 \Delta^{-1/2}(h) e^{-\lambda \rho(h)} 
\]
for 
$n\in\Ni$,  
$\lambda\in\Ri$
and
$g,h\in G$. 
Then, setting 
$h=e$, 
the required Gaussian estimate of 
$K^{(n)}$ 
follows by optimizing over 
$\lambda$ in the following 
bounds: 
\begin{eqnarray} 
e^{\lambda(\rho(g)-1)} \Delta^{1/2}(g) K^{(n)}(g) 
& \leq &  
\| \widetilde{T}_{\lambda}^n \|_{\widehat{1}\to \infty} 
  \nonumber \\
 & \leq & \| \widetilde{T}_{\lambda} \|_{\widehat{2}\to \infty} 
   \| \widetilde{T}_{\lambda}^{n-2} \|_{\widehat{2}\to \widehat{2}} 
   \| \widetilde{T}_{\lambda}\|_{\widehat{1}\to \widehat{2}} 
\leq c e^{\omega \lambda^2 n}
\label{ettilde} 
\end{eqnarray} 
for all 
$g\in G$,   
$n\geq 3$ 
and 
$\lambda\in\Ri$. 
Here,  
$\|\cdot\|_{\widehat{p}\to \widehat{q}}$ 
denotes the norm of an operator from 
$L^p(G;d\widehat{g})$ to 
$L^q(G;d\widehat{g})$.  
To justify (\ref{ettilde}), 
note that 
estimates of type  
$\| \widetilde{T}_{\lambda}^n\|_{\widehat{2}\to \widehat{2}} 
\leq 
e^{\omega \lambda^2 n}$, 
$n\in\Ni$,  
are a consequence of Theorem~\ref{tcentgauss1}  
since $\rho\in {\cal E}$ and 
$\Delta^{1/2}\colon L^2(G;dg) \to L^2(G;d\widehat{g})$ 
is unitary.  
The bounds of type  
$\| \widetilde{T}_{\lambda} \|_{\widehat{2}\to \infty}
+  \| \widetilde{T}_{\lambda} \|_{\widehat{1}\to \widehat{2}}
\leq c e^{\omega \lambda^2}$, 
$\lambda\in\Ri$, 
follow by integration of (\ref{ekgauss1}).  
We omit further details.

Theorem~\ref{tcentdiff1} is proved in the same way as  
\cite[Theorem~1.11]{Dun20}.  
In this proof, 
the estimate 
$\|(I-T_{\lambda}) T_{\lambda}^n \|_{2\to 2}
\leq c n^{-1} e^{\omega \lambda^2 n}$ 
(where 
$T_{\lambda}: = e^{\lambda\psi} T e^{-\lambda\psi}$)  
is deduced from the analyticity of 
$T$,  
by applying a perturbation theorem for analytic operators of 
Blunck \cite{Blu1}.  
For this deduction one needs the perturbation estimate~(\ref{epert3}). 

Alternatively, a simpler proof of the first estimate of 
Theorem~\ref{tcentdiff1} 
follows from the technique of \cite{Dun30}.  
In fact, 
given Theorems~\ref{tcentanal1} and \ref{tcentgset} the 
argument of \cite{Dun30} yields 
\[
\| \chi_E (I-T) T^n \chi_F \|_{2\to 2} 
\leq cn^{-1} e^{-b d(E,F)^2/n}  
\]
for all 
$n\in\Ni$  
and all Borel sets 
$E,F\subseteq G$.  
By Proposition~\ref{pgexp},    
this bound is equivalent 
to the desired estimate of 
$\|(I-T_{\lambda}) T_{\lambda}^n \|_{2\to 2}$.  

To obtain the remaining estimates of Theorem~\ref{tcentdiff1}, 
one observes that 
\begin{eqnarray*} 
\| \partial_h f\|_2^2 
& \leq & \rho(h)^2 \Gamma_2(f)^2     \\
& \leq & c' \rho(h)^2 \RRe ((I-T)f,f)  \\
& \leq & c \rho(h)^2
 \left[ \RRe ((I-T_{\lambda})f,f) + c \lambda^2 \|f\|_2^2 \right] 
\end{eqnarray*} 
for all $h\in G$, $f\in L^2$ and $|\lambda|\leq 1$, 
where the last two steps follow from 
Lemma~\ref{lrealest}
together with (\ref{epert3}).  
After replacing $f$ with $T_{\lambda}^n f$ 
and applying the first estimate of Theorem~\ref{tcentdiff1}, 
this enables one to get the desired estimate 
of
$\| \partial_h T_{\lambda}^n \|_{2\to 2}$.  
The similar estimate on 
\[
\| e^{\lambda\psi} \partial_h T^n e^{-\lambda\psi} \|_{2\to 2} 
= \| (e^{\lambda\psi} \partial_h e^{-\lambda\psi}) 
T_{\lambda}^n \|_{2\to 2} 
\]
can then be deduced by a straightforward argument, 
via the operator identity 
$e^{\lambda\psi} \partial_h e^{-\lambda\psi} 
= \partial_h + [ e^{-\lambda \partial_h \psi} -1] L(h)$.

As for Theorem~\ref{tpolygauss}, 
the required Gaussian bounds for  
$K^{(n)}(g)$ 
and 
$\partial_h K^{(n)}(g)$ 
follow from 
(\ref{ekgauss1}), 
Theorem~\ref{tcentgauss1} and
Theorem~\ref{tcentdiff1}, 
by choosing 
$\psi=\rho\in {\cal E}$    
and  applying 
\cite[Theorem~2.3]{Dun10}.

Integration of the Gaussian bound on $K^{(n)}$ 
in Theorem~\ref{tpolygauss}  
implies an estimate of type 
\[
\| e^{\lambda\rho} K^{(n)} \|_2  \leq c n^{-D/4} e^{\omega \lambda^2 n} 
\]
for all 
$n\in\Ni$ and 
$\lambda\geq 0$.   
Then, 
one may deduce 
the desired Gaussian bound of Theorem~\ref{tpolygauss} 
for  
$K^{(n)}-K^{(n+1)}$,  
by using the first estimate of Theorem~\ref{tcentdiff1} 
together with a factorization 
\begin{eqnarray*} 
& & e^{\lambda\rho(g)} | K^{(3m)}(g) - K^{(3m+1)}(g)|  \\
& \leq & \| e^{\lambda\rho} (K^{(2m)} - K^{(2m+1)}) \|_2 
         \| e^{\lambda\rho} K^{(m)} \|_2      \\
& \leq & 
 \| e^{\lambda \rho} (I-T)T^m e^{-\lambda\rho} \|_{2\to 2} 
 \| e^{\lambda \rho} K^{(m)} \|_2 \| e^{\lambda \rho} K^{(m)} \|_2 
\end{eqnarray*} 
for 
$g\in G$, 
$m\in\Ni$ and 
$\lambda\geq 0$.  
Here, we applied the identities 
$K^{(3m)}-K^{(3m+1)} = (K^{(2m)} - K^{(2m+1)}) * K^{(m)}$ 
and 
$K^{(2m)}-K^{(2m+1)} = (I-T)T^m (K^{(m)})$. 
Then optimize over 
$\lambda$. 
(For similar arguments,   
see \cite{Blu1}   
and the proof of 
\cite[Theorem~1.12]{Dun20}.)

To prove Theorem~\ref{tunifgauss},  
let the centered densities 
$(K_n)_{n=1}^{\infty}$ 
and the associated Markov operators 
$T_n f: = K_n *f$  
be as in the hypothesis.  
Since the $K_n$ satisfy the assumptions uniformly in $n$, 
an inspection of the proofs of Sections~\ref{s3} and \ref{s4} 
shows that the estimate 
\begin{equation} 
\| e^{\lambda\psi} T_n e^{-\lambda\psi}\|_{2\to 2} 
\leq e^{\omega \lambda^2} 
\label{euniftn} 
\end{equation} 
holds uniformly for all 
$n\in\Ni$, $\lambda\in\Ri$ and 
$\psi\in {\cal E}$.   
To see this, 
one observes that (\ref{euppgradpsi}) and  
the estimates of Lemma~\ref{lcrudest} and Proposition~\ref{ppert} 
hold with 
$T_n$ replacing $T$, 
with constants independent of $n$.  

Now (\ref{eunifg1})  is an immediate consequence of 
(\ref{euniftn}).  
The Gaussian bound (\ref{eunifg2}) then follows by the same 
technique used to prove Theorem~\ref{tcentgauss2}.

We next outline the proof of Theorem~\ref{tpolyunif}. 
The uniform estimate 
\[
\|K_{m+n} * \cdots * K_{m+1}\|_{\infty} 
\leq c n^{-D/2} 
\] 
for all 
$m\in\Ni_0$, $n\in\Ni$,   
follows from the proof of  
\cite[Theorem~VII.1.2]{VSC}.    
Together with the bounds~(\ref{eunifg1}), 
this enables one to run the argument 
of \cite[Lemma~2.2]{HebSal} to derive 
the Gaussian bound of Theorem~\ref{tpolyunif}.   
We omit further details.

Finally, let us prove Theorem~\ref{tunifanal}. 
With notation as in the theorem, 
consider the ``numerical range''  
\[
\Theta_{\alpha}: = \{ (T_{(\alpha)}f,f)\colon 
f\in L^2, \|f\|_2=1\} \subseteq \Ci 
\] 
of $T_{(\alpha)}$ 
for each 
$\alpha\in A$. 
Because of the uniform assumptions on the 
$K_{(\alpha)}$,  
an inspection of the proofs of 
Sections~\ref{s2} and \ref{s3} 
shows that estimate~(\ref{esect1}) 
holds for 
$T_{(\alpha)}$ 
with a constant 
$c>0$ independent of $\alpha$. 
This means that for some  
constant 
$\theta\in (0,\pi/2)$ 
one has 
\[
\Theta_{\alpha} \subseteq 
  \{\lambda\in \Ci\colon |\arg(1-\lambda)|
 \leq \theta\} 
\]
for all 
$\alpha\in A$.  
Similarly, 
the proof of Proposition~\ref{pspect} 
yields for some 
$\delta\in (0,1]$ 
that 
$\Theta_{\alpha}\subseteq \Lambda_{\delta}$ 
for all 
$\alpha\in A$,    
where 
$\Lambda_{\delta}$ is defined as in (\ref{elambdadef}). 
Hence 
\[
\Theta_{\alpha} \subseteq \Gamma_{\delta,\theta}: 
= \Lambda_{\delta} \cap \{\lambda \in \Ci\colon 
  |\arg(1-\lambda)| \leq \theta\},   
\]
and from Lemma~\ref{lspectr} we infer a resolvent estimate 
\begin{equation} 
\| (\lambda I- T_{(\alpha)})^{-1} \|_{2\to 2} \leq (dist(\lambda, 
\Gamma_{\delta,\theta}))^{-1} 
\label{eunifres} 
\end{equation} 
uniformly for all 
$\lambda\in \Ci\backslash \Gamma_{\delta,\theta}$  
and 
$\alpha\in A$.  
Using (\ref{eunifres}), 
one may check that the argument of 
\cite[p.102]{N1}
applies uniformly in 
$\alpha$, 
to give 
$\|(I-T_{(\alpha)}) T_{(\alpha)}^n \|_{2\to 2} \leq c n^{-1}$ 
for all 
$\alpha\in A$ and 
$n\in\Ni$.  
(The argument of \cite{N1} is essentially the proof 
of the implication 
(III)$\Rightarrow$(I) 
of Theorem~\ref{tanalcrit}.)  
This establishes Theorem~\ref{tunifanal} when 
$p=2$. 

Then, 
thanks to (\ref{eunifres}), 
one can verify that the proof of the interpolation 
theorem for analytic operators (\cite[Theorem~1.1]{Blu2}) 
goes through uniformly in 
$\alpha$, 
to yield Theorem~\ref{tunifanal} for 
$p\in (1,\infty)$.  
We omit further details.  
\hfill $\Box$

\section{Sublaplacians}  \label{s8}

In this section, 
we prove 
Theorems~\ref{tsublaps} 
and \ref{tsubpoly}.  

\bigskip

\noindent {\bf Proof of Theorem~\ref{tsublaps}.} 
Let 
$H=-\sum_{i=1}^{d'}  A_i^2 + A_0$ 
be a sublaplacian with drift on the connected Lie group 
$G$, 
and write 
$S_t: = e^{-tH}$ 
and 
$S_t f = K_t *f$ 
for  
$t>0$. 
It is well known that  
$(t,g)\mapsto K_t(g)$  
is a strictly positive, 
$C^{\infty}$ function on  
$(0,\infty)\times G$. 
Moreover, 
$K_t$ is a sub-Gaussian density on $G$ for each fixed 
$t>0$ 
(see,  for example,  
\cite[Appendix~A.4]{Var6}), 
and one has the      
following estimates.  
Given a compact interval 
$[\delta_1, \delta_2] \subseteq (0,\infty)$, 
$k, m\in\Ni_0 = \{0,1,2, \ldots\}$
and any right invariant vector fields 
$X_1, \ldots, X_m$,  
then there exist 
$c,b>0$ 
with 
\begin{equation} 
|(d/dt)^k X_1 \ldots X_m K_t(g)| 
\leq c e^{-b\rho(g)^2} 
\label{egaussloc} 
\end{equation} 
for all 
$t\in [\delta_1,\delta_2]$ 
and 
$g\in G$. 
In fact, 
(\ref{egaussloc}) follows from the 
sub-Gaussian property 
by applying a local 
parabolic Harnack inequality of the type given in  
\cite[Theorem~III.2.1]{VSC}. 

Assuming that 
$H$ is centered,  
one easily checks that
$A_0 \pi_1^{(j)}=0$ 
and 
$H \pi_1^{(j)}=0$ 
for all 
$j\in \{1, \ldots, q\}$.  
It is straightforward to deduce that 
$K_t * \pi_1^{(j)} = \pi_1^{(j)}$,  
which implies that 
$K_t$ is a centered density for each 
$t>0$.  

These observations imply that 
the family of densities 
$\{K_u \colon 1 \leq u \leq 2\}$ 
satisfies the hypotheses of Theorem~\ref{tunifanal}.  
Fix 
$p\in (1,\infty)$.  
Theorem~\ref{tunifanal} gives an estimate  
\[
\| (I-S_u) (S_u)^n \|_{p\to p} 
= \| (I-S_u) S_{nu}\|_{p\to p} \leq cn^{-1} 
\] 
for all 
$n\in\Ni$ and 
$u\in [1,2]$.   
Since 
$\|S_t\|_{p\to p} \leq 1$ for all $t$,  
it is easy to deduce that 
\begin{equation} 
\|(I-S_s) S_t \|_{p\to p} \leq c' t^{-1}  
\label{esdiff} 
\end{equation} 
for all 
$t\geq 2^{-1}$ and 
$s\in [0,1]$.  
Note that (\ref{egaussloc}) implies that 
$H S_s = - (d/ds) S_s$ 
is bounded in 
$L^p$ for every 
$s>0$. 
By writing 
$I-S_1 = \int_0^1 ds\, HS_s$ 
and 
\[
H S_t 
 = (I-S_1)S_t + HS_{1/2} \int_0^1 ds\, (I-S_s) S_{t-(1/2)}, 
\]
we obtain 
\[
\|H S_t\|_{p\to p} 
\leq ct^{-1} + c' \| HS_{1/2}\|_{p\to p} t^{-1} 
\leq c'' t^{-1}  
\]
for all 
$t\geq 1$.  
Theorem~\ref{tsublaps} is proved. 
\hfill $\Box$ 

\bigskip 

\noindent {\bf Remark.} 
The above proof establishes 
the following abstract result.  
Let 
$S_t = e^{-tH}$, $t\geq 0$, 
be a uniformly bounded semigroup 
of operators in a Banach space. 
Suppose (\ref{esdiff}) 
holds for all 
$s\in [0,1]$ and 
$t\geq t_0$,  
for some constant 
$t_0>0$, 
and that  
$HS_u$ is bounded for some $u>0$. 
Then
$\| HS_t\|\leq ct^{-1}$ 
for all sufficiently large $t$. 

It is not clear to us whether 
the bound  
$\|HS_t\|\leq ct^{-1}$ 
(for large $t$) 
could be derived 
from a weaker version 
$\|(I-S_1)S_t\| \leq ct^{-1}$, $t\geq 1$, 
of (\ref{esdiff}).  
If this were possible then    
we would only need Theorem~\ref{tcentanal1}, 
and not 
Theorem~\ref{tunifanal}, 
to prove Theorem~\ref{tsublaps}. 

\bigskip

\noindent {\bf Proof of Theorem~\ref{tsubpoly}.}  
Theorem~\ref{tpolygauss} applied to the density 
$K=K_1$ 
gives 
\begin{equation} 
K_t(g) + t^{1/2} |(\partial_h K_t)(g)| 
\leq c t^{-D/2} e^{-b\rho(g)^2/t} 
\label{eroughbd} 
\end{equation} 
for all 
$g\in G$, $h\in U$ 
and 
$t\in \Ni= \{1,2,3, \ldots \}$.    
But by the convolution 
identities 
$K_{t_1+t_2}= K_{t_1}* K_{t_2}$, 
$\partial_h K_{t_1+t_2}
= (\partial_h K_{t_1})* K_{t_2}$ 
and by 
(\ref{egaussloc}), 
the bound (\ref{eroughbd}) easily extends to all 
real 
$t\geq 1$. 

Next, 
applying Theorem~\ref{tsublaps} and (\ref{eroughbd}) we have 
\begin{eqnarray} 
\| (d/dt)^2 K_t\|_{\infty} & = & \| H^2 K_t\|_{\infty}  
= \| (HK_{t/2}) * (HK_{t/2})\|_{\infty} \nonumber   \\
& \leq & \| HK_{t/2}\|_2^2   \nonumber \\
& \leq & \| HS_{t/4}\|_{2\to 2}^2 \| K_{t/4}\|_2^2 
\leq c t^{-2} t^{-D/2} 
 \label{eunifbd1} 
\end{eqnarray} 
for all 
$t\geq 4$.   
A Taylor expansion of the function 
$t\mapsto K_t(g)$ 
shows that
\[
| (d/dt) K_t(g)| 
\leq s^{-1} |K_{t+s}(g)| + s^{-1} |K_t(g)| 
+ 2^{-1} s \sup_{u\in [t,t+s]} |(d/du)^2 K_u(g)| 
\]
for all 
$t,s>0$ and 
$g\in G$. 
By choosing 
$s= t e^{-\varepsilon \rho(g)^2/t}$ 
for a sufficiently small constant 
$\varepsilon>0$, 
and applying (\ref{eroughbd}), (\ref{eunifbd1}), 
we deduce that 
\begin{equation} 
|(d/dt)K_t(g)| \leq c' t^{-1} t^{-D/2} e^{-b' \rho(g)^2/t} 
\label{etimebd1} 
\end{equation} 
for 
$t\geq 4$.   

Finally, 
let $X$ be a right invariant vector field on $G$. 
A standard local regularity estimate for subelliptic operators 
(cf. \cite[Corollary~III.1.3]{VSC})
implies that  
\[
|X F(0,e)| 
\leq c \sup_{(s,h)\in (-1,1)\times U} 
  |F(s,h)| 
\]
for all solutions 
$F\colon (-1,1)\times U \to \Ri$ 
of the heat equation 
$(\partial/(\partial t) + H)F=0$.   
Applying this to the functions 
$F^{(t,g)}$ defined by 
\[
F^{(t,g)}(s,h): = K_{t+s}(hg)-K_t(g) 
\]
gives 
\[
|XK_t(g)| \leq c \sup_{(s,h)\in (-1,1)\times U} 
   |K_{t+s}(hg) - K_t(g)| 
\]
for all 
$t\geq 2$ and 
$g\in G$. 
The desired Gaussian bound on $XK_t(g)$ 
then follows easily from (\ref{eroughbd}) and (\ref{etimebd1}). 
\hfill $\Box$

\section{Locally compact groups} \label{s9} 

In this section and the next, 
our aim is to prove the results of Subsection~\ref{s1s3}.   
In general, 
$G$ will denote a compactly generated locally compact group, 
and 
$G_0 = \overline{[G,G]}$, $G_1$ 
are  
the closed normal subgroups of 
$G$ defined as in 
Subsection~\ref{s1s3}.

We fix some notation.   
For any locally compact group $H$, 
the connected component containing the identity of $H$ 
will be denoted $H^c$.  
Note that $H^c$ is always a closed normal subgroup of $H$ 
and the quotient $H/H^c$ is a totally disconnected 
locally compact group 
(cf. \cite[Chapter~II]{HR}); 
$H/H^c$ is discrete if and only if $H^c$ is open in $H$, 
which is not always the case.  
Recall that $H$ is said to be {\it almost connected} 
if $H/H^c$ is compact.   

Lie groups are {\it not} assumed to be connected in what follows. 
If $H$ is a Lie group, however,  
then $H^c$ is an open subgroup of $H$, 
and $H/H^c$ is a countable discrete group.

We shall need the following fundamental structure theorems for locally 
compact groups, 
for which \cite{MZ} is a standard reference.

\begin{thm} \label{tmz1}
(\cite[Section~4.6]{MZ}.)     
Let $H$ be an almost connected locally compact group. 
Then there exists a compact normal subgroup $N$ of $H$ 
such that 
$H/N$ is a Lie group. 
\end{thm} 

\begin{thm} \label{tmz2} 
(\cite[Section~2.3]{MZ}.)   
Let $H$ be a locally compact group. 
There exists an open subgroup 
$H'$ of $H$ such that 
$H'$ is almost connected.  
\end{thm} 

By a {\it one-parameter subgroup} of a locally compact group 
$H$ we mean a continuous homomorphism 
$\theta\colon \Ri\to H$.

\begin{thm} \label{tmz3} 
(\cite[Section~4.15]{MZ}.)    
Let $H$ be a locally compact group, 
$N$ a compact normal subgroup of $H$, 
and  
$p\colon H\to H/N$ 
the canonical homomorphism.   
If 
$\widetilde{\theta}$ is a one-parameter subgroup of 
$H/N$,  
then there exists a one-parameter subgroup 
$\theta$
of $H$ such that 
$\widetilde{\theta}(t) = p(\theta(t))$ 
for all $t\in\Ri$. 
\end{thm}

\noindent {\bf Proof of Theorem~\ref{tlcanal}.} 
Let $K$ be centered and satisfy Assumption~(i)'.  
Just as in the proof of Lemma~\ref{lrealest} 
(see the Remark following Lemma~\ref{lgradbas}) 
we have an estimate 
\[
\RRe ((I-T)f,f) \geq c \Gamma_2(f)^2
\]
for all 
$f\in L^2$. 
Then analyticity of $T$ in $L^2$ 
is a consequence of the arguments of Section~\ref{s2} together 
with (\ref{eupgd22}). 
Analyticity in $L^p$ ($1<p<\infty$) then follows as in 
Corollary~\ref{ccentanalp}. 

Therefore, 
it only remains to prove (\ref{eupgd22}).   
For this we need the following version of (\ref{eg1gen}).   

\begin{prop} \label{pg0g1ext} 
Let $C$ be any compact subset of $G_1$. 
Then there exists 
$c=c(C)>0$ such that 
\[
|(\partial_g f_1,f_2)| \leq c \Gamma_2(f_1) \Gamma_2(f_2)  
\]
for all $g\in C$ and $f_1,f_2\in L^2$. 
\end{prop} 

The proof of Proposition~\ref{pg0g1ext}
is rather technical and is deferred to Section~\ref{s10}.  

Next, 
consider the homomorphism 
$\pi_1\colon G\to G/G_1 \cong \Ri^q \times \Zi^r$, 
and identify 
$G/G_1= 
\Ri^q\times \Zi^r$. 
Let 
$e_1, \ldots, e_{q+r}\in \Zi^{q+r} 
\subseteq \Ri^q\times \Zi^r$
be  the standard basis vectors,   
that is, 
$e_j$ has a $1$ in the $j$-th position and zeroes elswhere. 

As in the proof of Theorem~\ref{tgradupp1} 
for Lie groups, 
we wish to find one-parameter subgroups with the following 
property.

\begin{lemma} \label{lchooseone} 
There exist one-parameter subgroups $x_1, \ldots, x_q$ of $G$ 
such that 
$\pi_1(x_j(t)) = t e_j$ 
for all 
$t\in\Ri$ 
and 
$j\in \{1, \ldots, q\}$.  
\end{lemma} 
\proof\ 
By general topological group theory,  
$\pi_1(G^c)$ 
is a dense, connected subgroup of 
the group 
$(G/G_1)^c = \Ri^q\times \{0\}$  
(see \cite[Theorem~II.7.12]{HR}).  
Applying Theorem~\ref{tmz1} to 
$G^c$, 
we choose a compact normal subgroup $N$ of $G^c$ 
such that 
$H': = G^c/N$
is a Lie group.
Note that 
$N\subseteq G_1$ by compactness of $N$.   
Therefore, 
the homomorphism 
$\pi_1|_{G^c}\colon G^c\to \Ri^q \times \{0\}$ 
induces a homomorphism 
$\pi_2\colon H' \to \Ri^q\times \{0\}$ 
with 
$\pi_2(gN) = \pi_1(g)$, 
$g\in G^c$.  
Thus 
$\pi_1(G^c)= \pi_2(H')$ is a dense, connected 
Lie subgroup of 
$\Ri^q\times \{0\}$, 
and hence 
$\pi_1(G^c) = \pi_2(H') = \Ri^q\times \{0\}$. 

Since $H'$ is a Lie group, 
as in the proof of Theorem~\ref{tgradupp1} 
one can find one-parameter subgroups 
$y_1, \ldots, y_q$ 
of $H'$ 
such that  
$\pi_2(y_j(t))= te_j$,  $j\in \{1,\ldots, q\}$.  
Applying Theorem~\ref{tmz3} to  
$G^c$ 
yields the existence of one-parameter subgroups 
$x_1, \ldots, x_q$ of $G^c$ such that 
$y_j(t) = x_j(t) N$ for all $t\in\Ri$.  
These subgroups satisfy the lemma. 
\hfill $\Box$ 

\bigskip 

Next, 
fix elements 
$z_1, \ldots, z_r\in G$ 
such that 
$\pi_1(z_j) = e_{q+j}$
for 
$j\in \{1, \ldots, r\}$. 
Given any 
$g\in G$, we set 
$t_j = \pi_1^{(j)}(g)$,  $j\in \{1, \ldots, q+r\}$, 
and write 
\[
g = x_1(t_1) \cdots x_q(t_q) 
    z_1^{t_{q+1}} \cdots z_r^{t_{q+r}} 
    g' 
\]
where $g'\in G$.  
Then 
$g'\in G_1$, because 
\[
\pi_1(g) = (t_1, \ldots, t_{q+r}) 
= \pi_1(x_1(t_1) \cdots x_q(t_q) z_1^{t_{q+1}} \cdots z_r^{t_{q+r}}). 
\]
Let 
$C\subseteq G$ be any compact set.    
By applying identity~(\ref{ediffid2}) and 
arguing as in the proof of Theorem~\ref{tgradupp1}, 
one finds that 
\[
\partial_g = \sum_{j=1}^q \pi_1^{(j)}(g) \partial_{x_j(1)} 
           + \sum_{k=1}^r \pi_1^{(q+k)}(g) \partial_{z_k} 
           + V_g,  
\]
where $V_g$ satisfies an estimate 
\begin{equation} 
|(V_g f_1,f_2)| \leq c(C)  \Gamma_2(f_1) \Gamma_2(f_2) 
\label{elcvest} 
\end{equation} 
for all $f_1,f_2\in L^2$ and all 
$g\in C$, 
with $c(C)>0$ a constant depending on $C$. 
The proof of (\ref{elcvest}) follows the argument 
of Theorem~\ref{tgradupp1}  
with the following main differences:  
first, 
one now applies 
Proposition~\ref{pg0g1ext} to bound the term  
$(\partial_{g'} f_1,f_2)$,  
and secondly, 
one requires estimates of form  
\[
|((\partial_{z_j^m}- m \partial_{z_j})f_1,f_2)|  
\leq c(m) \Gamma_2(f_1)\Gamma_2(f_2) 
\]
for all 
$j\in \{1, \ldots, r\}$, 
$m\in \Zi$ and 
$f_1,f_2 \in L^2$. 
The latter estimates are easily obtained,  
since using (\ref{ediffid2}) and the identity 
$\partial_{z_j^{-1}} 
= - \partial_{z_j} - \partial_{z_j} \partial_{z_j^{-1}}$ 
one may express 
$\partial_{z_j^m} - m \partial_{z_j}$ 
as a finite linear combination of terms each of form 
\[
\partial_{z_j^{\varepsilon_1}} \ldots \partial_{z_j^{\varepsilon_l}}
\]
where 
$l\geq 2$ and 
$\varepsilon_1, \ldots \varepsilon_l \in \{-1,+1\}$.  

Finally, since $K$ is centered, 
by choosing $C$ to contain the compact support of $K$ then 
(\ref{elcvest}) yields 
\[
|((I-T)f_1,f_2)| =  
\left|\int_G dg\, K(g) (V(g)f_1,f_2) \right| 
\leq c \Gamma_2(f_1)\Gamma_2(f_2),  
\]
proving (\ref{eupgd22}).    
This ends the proof of Theorem~\ref{tlcanal}, 
modulo Proposition~\ref{pg0g1ext}. 
\hfill $\Box$

\bigskip

\noindent {\bf Proof of Theorem~\ref{tlcgauss}.} 
It is enough to prove Statement~(I) of the theorem, that is, 
the estimate~(\ref{eupgdinf1}).   
For then Statements~(II) and (III) can be deduced just as in 
the proofs of Theorems~\ref{tcentgauss1}, 
\ref{tcentgset}  
and \ref{tcentgauss2}.   

To obtain (\ref{eupgdinf1}),   
we need the following variation of   
Proposition~\ref{pg0g1ext}. 

\begin{prop} \label{pg0g1ext2} 
Let $C\subseteq G_1$ be compact. 
Then there exists $c=c(C)>0$ such that 
\[
|(\partial_g f_1,f_2)| \leq c \Gamma_{\infty}(f_1) \Gamma_1(f_2) 
\]
for all $g\in C$ and
$f_1\in L^{\infty}$, 
$f_2\in L^1$.  
\end{prop} 

Reasoning as in the proof of Theorem~\ref{tlcanal} then shows, 
analogously to (\ref{elcvest}), 
that given a compact 
$C\subseteq G$ one has 
\[
|(V_g f_1,f_2)| \leq c(C) \Gamma_{\infty}(f_1) \Gamma_1(f_2) 
\]
for all 
$g\in C$, 
$f_1\in L^\infty$, 
$f_2\in L^1$. 
Then (\ref{eupgdinf1}) follows using the hypothesis that 
$K$ is centered. 
Theorem~\ref{tlcgauss} is proved.  
\hfill $\Box$

\bigskip

\noindent {\bf Proof of Theorem~\ref{tlcpolygauss}.}   
This follows the proof of Theorem~\ref{tpolygauss}. 
\hfill $\Box$  

\bigskip

\noindent {\bf Proof of Theorem~\ref{tamengen}.}  
The proof of Theorem~\ref{tamengen} follows that of Theorem~\ref{tamen}. 
The only significant difference occurs 
when proving that Conditions~(II) or (IV) fail for non-centered $K$; 
one must also consider the case where $P$ 
(the image of $K$ under $\pi_1^{(j)}$) 
is a non-centered density 
on $\Zi$. 
For such $P$,  
the relevant analysis 
is already given  
in, for example, \cite{Dun20}. 
\hfill $\Box$ 

\bigskip

\noindent {\bf Proof of Theorem~\ref{tlcunif}.} 
The proof of Theorem~\ref{tlcgauss} 
goes through uniformly for the densities $(K_n)_{n=1}^{\infty}$, 
and the desired estimates follow 
just as in the proof of Theorem~\ref{tunifgauss}. 
\hfill $\Box$

\section{Proof of Propositions~\ref{pg0g1ext} and \ref{pg0g1ext2}} \label{s10} 

Let $G$ be locally compact and compactly generated.   
In this section, 
we will prove Proposition~\ref{pg0g1ext};   
the proof of Proposition~\ref{pg0g1ext2} requires only obvious 
changes and is omitted.

The proof of Proposition~\ref{pg0g1ext} reduces 
to the case where $G$ is second countable, 
by the following theorem and lemma.  

\begin{thm} \label{ttwocount} 
(\cite[Theorem~II.8.7]{HR}.)  
Let $H$ be a locally compact group which is 
$\sigma$-compact 
(that is, $H$ is the union of some countable family of 
compact subsets of $H$; 
in particular, any compactly generated group is 
$\sigma$-compact).    
Then there exists a compact normal subgroup $N$ of $H$ 
such that 
$H/N$ 
is a second countable locally compact group.  
\end{thm}

\begin{lemma} \label{ltwocred}
Let $N$ be a compact normal subgroup of $G$.   
If Proposition~\ref{pg0g1ext} holds with respect to the 
group $G': = G/N$, 
then it also holds with respect to 
$G$. 
\end{lemma} 

The proof of Lemma~\ref{ltwocred} is not difficult and is left to 
the reader.   

Thus, in the rest of this section we shall assume that 
$G$ is second 
countable.

The next step in the proof   
is to obtain the desired estimate in case  
$g\in G_0$, as follows.

\begin{prop} \label{plcg0} 
For each $g\in G_0$, 
there exists a $c(g)>0$ such that 
\begin{equation} 
|(\partial_g f_1,f_2)| \leq c(g) \Gamma_2(f_1) \Gamma_2(f_2) 
\label{egrup22g0} 
\end{equation} 
for all $f_1,f_2\in L^2$. 
\end{prop} 

The final step consists in extending the estimate from 
$G_0$ to $G_1$; 
this extension 
depends on the compactness of 
$G_1/G_0$.   

We begin the proof of Proposition~\ref{plcg0} with 
the following observation. 

\begin{lemma} \label{lggcp} 
Let $N$ be any compact subgroup of $G$. 
Then an estimate (\ref{egrup22g0}) holds for each $g\in [G,G]N
= N[G,G]$.  
\end{lemma} 
\proof\ 
This is an straightforward adaption of the proofs of 
Lemmas~\ref{luest} and \ref{luoneest}, 
using (\ref{ediffid}), (\ref{enabla1}) and 
(\ref{enabla2}).     
\hfill $\Box$ 

\bigskip 

Proposition~\ref{plcg0} would now follow if we had   
$G_0 \subseteq N[G,G]$ with $N$ some compact subgroup of $G$. 
Unfortunately we do not know if such an $N$ exists in general,  
but we have the following substitute result. 
If $A_1, \ldots, A_s$ are subsets of a group $H$, 
denote by 
$\langle A_1, \ldots, A_s\rangle$ 
the subgroup of 
$H$ generated by 
$A_1 \cup A_2 \cup \cdots \cup A_s$. 

\begin{prop} \label{pg0gg} 
There exist a compact subgroup $N$ of $G_0$, 
and one-parameter subgroups 
$\theta_1, \ldots, \theta_s$ of $G_0$ 
($s\in \Ni_0=\{0,1,2, \ldots\}$)  
such that 
$\theta_j(1) \in [G,G]N$
for all $j$ and 
\begin{equation} 
G_0 = [G,G]N \langle \theta_1(\Ri), \ldots \theta_s(\Ri) \rangle. 
\label{eg0rep} 
\end{equation}  
\end{prop} 

Proposition~\ref{pg0gg} 
is just a particular case of the following lemma 
(apply the lemma with 
$H=G_0=\overline{[G,G]}$ and 
$A=[G,G]$).    

\begin{lemma} \label{lhdense} 
Let $A$ be a dense subset of a locally compact group $H$. 
Then there exist a compact subgroup $N$ of $H$, 
and one-parameter subgroups 
$\theta_1, \ldots, \theta_s$ of 
$H$ 
($s\in\Ni_0$)   
such that 
\begin{equation} 
\theta_j(1) \in AN 
\label{etheta} 
\end{equation} 
for all $j$, 
the set 
\begin{equation} 
H': = N \langle \theta_1(\Ri), \ldots \theta_s(\Ri)\rangle 
\label{ehprime} 
\end{equation}
is an open subgroup of 
$H$ with $N$ a normal subgroup of 
$H'$,    
and
\begin{equation} 
H = AH' = AN \langle \theta_1(\Ri), \ldots, \theta_s(\Ri) \rangle. 
\label{ehrep} 
\end{equation} 
\end{lemma} 

\noindent {\bf Proof of Lemma~\ref{lhdense}.}  
Equation~(\ref{ehrep}) will follow from 
(\ref{ehprime}), 
because 
$H = \overline{A} =AV$ 
holds for any    
neighborhood $V$ of $e$ in $H$, 
hence for $V=H'$ since $H'$ is an open subgroup of $H$.

To prove the remaining statements of the lemma, 
we begin with some special cases of $H$.  
First,  suppose that 
$H$ is a Lie group (not necessarily connected)  
with Lie algebra 
$\gothh$.   
The exponential map
$\exp\colon \gothh\to H$ 
maps some open neighborhood of $0\in\gothh$ 
diffeomorphically onto a open neighborhood $W$ 
of $e$ in $H$. 
Since 
$A\cap W$ is dense in $W$, 
it is easy to deduce   
the existence of a vector space basis 
$X_1, \ldots, X_s$ of 
$\gothh$
such that 
$\exp(X_j)\in A\cap W$ 
for all $j$. 
Setting $\theta_j(t)= \exp(tX_j)$, $t\in\Ri$, 
we have 
$\theta_j(1)\in A$ 
and 
$H^c = \langle \theta_1(\Ri), \ldots, \theta_s(\Ri)\rangle$.  
Now 
$H^c$ is open in $H$ since $H$ is a Lie group, 
so 
(\ref{etheta}) and (\ref{ehprime}) hold 
with $N=\{e\}$ and $H' = H^c$. 

Next, suppose that $H$ is an almost connected locally compact group. 
By Theorem~\ref{tmz1}, choose a compact normal subgroup 
$N$ of $H$ 
such that 
$\widetilde{H}:= H/N$ 
is a Lie group, 
and let 
$p\colon H\to \widetilde{H}$ be the canonical map.  

Since $p(A)$ is dense in 
$\widetilde{H}$, 
by applying the preceding case we 
find one-parameter subgroups 
$\widetilde{\theta}_1, \ldots, \widetilde{\theta}_s$ of 
$\widetilde{H}$
such that 
$\widetilde{\theta}_j(1)\in p(A)$ 
and  
\[
\widetilde{H}^c = \langle \widetilde{\theta}_1(\Ri), 
      \ldots, \widetilde{\theta}_s(\Ri) \rangle. 
\]
By Theorem~\ref{tmz3} there exist one-parameter subgroups 
$\theta_1, \ldots, \theta_s$ of $H$ such that 
$\widetilde{\theta}_j(t) = p(\theta_j(t))$, 
$t\in\Ri$. 
Then 
$\theta_j(1) \in p^{-1}(p(A)) = AN$ for all 
$j$. 
Setting $H' := p^{-1}(\widetilde{H}^c)$, 
we observe that 
$H'$ is an open subgroup of $H$ 
and that 
$H' = N \langle \theta_1(\Ri), \ldots, \theta_s(\Ri)\rangle$.  
This proves the lemma for $H$ almost connected. 

Finally, 
for $H$ any locally compact group we choose by Theorem~\ref{tmz2} 
an open, almost connected subgroup 
$\widehat{H}$ of $H$. 
Since 
$\widehat{A}: = A\cap \widehat{H}$ is dense in 
$\widehat{H}$, 
we may apply the preceding case to 
$\widehat{H}$ 
to reach the desired conclusion. 
\hfill $\Box$ 

\bigskip 

\noindent {\bf Proof of Proposition~\ref{plcg0}.}  
Choose $N$ and 
$\theta_1, \ldots, \theta_s$ 
as in Proposition~\ref{pg0gg}.   
By Lemma~\ref{lggcp}, 
the desired estimate 
(\ref{egrup22g0}) 
holds when 
$g= \theta_j(1)$ for some 
$j$.  
One proves in the same way as 
(\ref{eauxe3}) that,  
for each $t\in\Ri$,  
there exists 
$c(t)>0$ with 
\[
|((\partial_{\theta_j(t)} - t \partial_{\theta_j(1)})f_1,f_2)| 
\leq c(t)   \Gamma_2(f_1)\Gamma_2(f_2) 
\]
for all 
$f_1, f_2\in L^2$.  
Therefore, 
(\ref{egrup22g0}) holds whenever 
$g\in \theta_j(\Ri)$.  
Then from (\ref{eg0rep}), 
from Lemma~\ref{lggcp} 
and by applying 
(\ref{ediffid}), (\ref{enabla1}), (\ref{enabla2}), 
it is straightforward to deduce (\ref{egrup22g0}) 
for each 
$g\in G_0$.  
\hfill $\Box$

\bigskip

We next show that the estimate of Proposition~\ref{plcg0}
self-improves into a uniform estimate over compacta.

\begin{lemma} \label{lg0cpct} 
Let $C$ be any compact subset of $G_0$. 
Then there exists 
$c=c(C)>0$ such that 
\[
|(\partial_g f_1,f_2)| \leq c \Gamma_2(f_1) \Gamma_2(f_2)  
\]
for all $g\in C$ and $f_1,f_2\in L^2$. 
\end{lemma} 
\proof\ 
For each $m\in\Ni$, 
let $C_m$ 
be the set of all $g\in G_0$ for which the estimate 
$ |(\partial_g f_1,f_2)| \leq m \Gamma_2(f_1) \Gamma_2(f_2)$ 
holds for all 
$f_1,f_2\in L^2$. 
By taking limits in this estimate, 
it is clear that  
$C_m$ is closed in $G_0$. 
Moreover,  
$G_0 = \bigcup_{m=1}^{\infty} C_m$ 
by Proposition~\ref{plcg0}.  
The Baire theorem for locally compact Hausdorff spaces 
(see \cite[Appendix~B.19]{HR}),   
applied to $G_0$, 
implies that 
$C_{m_0}$ has non-empty interior for some $m_0\in\Ni$. 

Select therefore a non-empty open subset $A$ of $G_0$ with compact closure 
$\overline{A}$, 
such that 
$\overline{A}\subseteq C_{m_0}$. 
Given a compact set 
$C\subseteq G_0$, 
we may find points 
$h_1, \ldots, h_k\in G_0$ such that 
$C \subseteq \bigcup_{i=1}^k h_i A$. 
For $g\in C$ write $g=h_i a$ where $a\in A$, 
and note that 
$\sup \{\rho(b)\colon b\in A\}  < \infty$ 
by compactness of $\overline{A}$.   
Then by Proposition~\ref{plcg0}, 
(\ref{ediffid}), (\ref{enabla1}) and (\ref{enabla2}) we obtain 
\[
|(\partial_g f_1, f_2)| 
\leq |(\partial_{h_i}f_1,f_2)| 
+ |(\partial_a f_1, L(h_i^{-1}) f_2)| 
\leq c \Gamma_2(f_1) \Gamma_2(f_2)      
\]               
where $c$ is independent of $g\in C$.  
\hfill $\Box$

\bigskip

To continue the proof of Proposition~\ref{pg0g1ext}, 
we exploit the compactness of the quotient group 
$G_1/G_0 \cong M$. 
Let us identify 
$M=G_1/G_0$.  
Since $G_0$ is a closed subgroup of the second countable, locally compact 
group $G_1$,  
by a standard result there exists a Borel set 
$\Sigma\subseteq G_1$ 
which meets each coset 
$g G_0$ 
($g\in G_1$)   
in exactly one point;   
since 
$M$ is compact,  we may choose 
$\Sigma$ to have compact closure in $G_1$ 
(see \cite[Lemma~1.1]{Mac1}).   

If we denote by 
$\pi\colon G_1\to G_1/G_0=M$ the canonical 
projection, 
then 
$\pi|_\Sigma\colon \Sigma\to M$ 
is a measurable bijection. 
Let 
$q\colon M \to \Sigma\subseteq G$ 
be the inverse of 
$\pi|_{\Sigma}$, 
so that 
$\pi(q(a)) = a$, 
$a\in M$. 
Consider the Haar measure $da$ on $M$, 
normalized so 
$da(M)=1$, 
and let 
$dy$ 
be the image measure of 
$da$ 
under 
$q\colon M\to \Sigma$.  
Then 
$dy(\Sigma)=1$,  
and 
\begin{equation} 
\int_{\Sigma} dy f(y) = \int_M da\, (f\circ q)(a) 
\label{eqrel} 
\end{equation} 
for all 
$f\in L^1(\Sigma; dy)$. 
The following result is essentially a generalization  
of the identity~(\ref{ecpctrep}).

\begin{lemma} \label{lcpsec} 
Let 
$g\in G_1$,  
and for each 
$y\in \Sigma$ 
define $z^y = z(g,y)\in G_1$ by 
\[
z^y : = q(\pi(yg^{-1})) (yg^{-1})^{-1}.    
\]
Then $z^y\in G_0$, and for each 
$f\in L^p=L^p(G;dg)$ 
one has 
\[
\partial_g f = - \int_{\Sigma} dy\, \partial_g \partial_{y^{-1}} f 
           + 
\int_{\Sigma} 
           dy\, L(q(\pi(yg^{-1}))^{-1}) \partial_{z^y} f. 
\]
\end{lemma} 
\proof\ 
The definition of 
$z^y$ implies  that 
$\pi(z^y)$ is the identity element 
of $M=G_1/G_0$, 
so that 
$z^y\in G_0$.  

It suffices to prove that the identity for $\partial_g f$ 
holds pointwise 
for any continuous function 
$f\colon G\to \Ci$. 
We start  by observing that  
\begin{eqnarray*} 
(\partial_g f)(e) & = & f(g^{-1})-f(e)    \\
 & = & f(g^{-1}) - \int_{\Sigma} dy\, f(yg^{-1})  
      - f(e) +   \int_{\Sigma} dy\, f(y)    \\
 &   & {} + \int_{\Sigma} dy\, f(yg^{-1}) - \int_{\Sigma} dy\, f(y) \\
 & = & - \int_{\Sigma} dy\, 
(\partial_g \partial_{y^{-1}}f)(e)     \\
& & {} +  \int_{\Sigma} dy\, [f(yg^{-1}) - f(q(\pi(yg^{-1})))]   \\ 
& & {} + \int_{\Sigma} dy\, f(q(\pi(yg^{-1}))) - \int_{\Sigma} dy\, f(y), 
\end{eqnarray*} 
where we used that  
$dy(\Sigma)=1$.  
The last two integrals in the right side cancel each other; 
indeed,  
with 
$F:= f\circ q$ we have from 
(\ref{eqrel}) that 
$\int_{\Sigma} dy\, f(y) = \int_M da \, F(a)$ 
and 
\[ 
\int_{\Sigma} dy\, f(q(\pi(yg^{-1}))) 
           = \int_M  da\, F(a \pi(g^{-1})) 
           = \int_M da\, F(a).  
\]
Because 
\[
\int_{\Sigma} dy\, [f(yg^{-1}) - f(q(\pi(yg^{-1})))]   
= \left( \int_{\Sigma} 
           dy\, L(q(\pi(yg^{-1}))^{-1}) \partial_{z^y} f \right)(e), 
\]
we obtain the desired identity at the point $e\in G$. 
The identity at any point $h\in G$ then follows by replacing 
$f$ with its right translate 
$R_G(h)f$ where 
$(R_G(h)f)(k) = f(kh)$, $k\in G$.  
\hfill $\Box$

\bigskip

Now suppose that $C\subseteq G_1$ is compact and $g\in C$. 
By Lemma~\ref{lcpsec}, 
\[
|(\partial_g f_1,f_2)| 
\leq 
\int_{\Sigma} dy\, | (\partial_{y^{-1}} f_1, \partial_{g^{-1}}f_2)| 
+ 
\int_{\Sigma} dy\, |(\partial_{z^y} f_1, 
                    L(q(\pi(yg^{-1}))) f_2)| 
\]
for $f_1,f_2\in L^2$. 
Then Proposition~\ref{pg0g1ext}  
follows by applying 
(\ref{enabla1}) to estimate the first integral over $\Sigma$, 
and 
applying Lemma~\ref{lg0cpct} and (\ref{enabla2}) 
to estimate the second integral since $z^y\in G_0$.  
Here we must observe that,  
because $\Sigma$ has compact closure 
and $y, q(\pi(yg^{-1}))\in \Sigma$, 
the elements 
$z^y= z(g,y)$ remain within a fixed compact 
subset of $G_0$ as $g$ varies over $C$. 
The proof of Proposition~\ref{pg0g1ext} is complete. 
\hfill $\Box$

\section*{Acknowledgements} 

During this work the author was located at the University of New 
South Wales, and was financially  supported 
by the Australian Research Council (ARC) 
Centre of Excellence for Mathematics and Statistics of Complex 
Systems (MASCOS).  

I thank George Willis for valuable advice on the structure theory 
of locally compact groups.

\bigskip 

\noindent Department of Mathematics   

\noindent Macquarie University 

\noindent NSW 2109, Australia

\noindent E-mail: ndungey@ics.mq.edu.au


\begin{thebibliography}{99}


\bibitem{Ale3} G. Alexopoulos, 
 `Sous-laplaciens et densit\'{e}s centr\'{e}s 
 sur les groupes de Lie \`{a} croissance polynomiale 
 du volume', 
 C. R. Acad. Sci. Paris S\'{e}r. I Math. 
 326 (1998),  539-542.  

\bibitem{Ale4} G. Alexopoulos, 
 `Sub-Laplacians with drift on Lie groups of polynomial volume 
 growth', 
 Mem. Amer. Math. Soc. 155 (2002), number 739.  

\bibitem{Ale5} G. Alexopoulos, 
 `Random walks on discrete groups of polynomial volume growth', 
 Ann. Probab. 30 (2002), 723-801.   

\bibitem{Ale6} G. Alexopoulos, 
 `Centered densities on Lie groups of polynomial volume growth', 
 Probab. Theory Relat. Fields 124 (2002), 112-150.  

\bibitem{Blu1} S. Blunck, 
 `Perturbation of analytic operators and temporal regularity 
  of discrete heat kernels', 
 Coll. Math. 86 (2000), 189-201.  

\bibitem{Blu2} S. Blunck,  
 `Analyticity and discrete maximal regularity
 on $L_p$-spaces', 
 J. Funct. Anal. 183 (2001), 211-230.  

\bibitem{Boun} N. Bounechada, 
 `Distorsion des distances dans les groupes de Lie nilpotents', 
 Bull. Sci. Math. 127 (2003), 797-813.  


\bibitem{Car} T. K. Carne, 
 `A transmutation formula for Markov chains', 
 Bull. Sci. Math. 109 (1985),  
 399-405.  




\bibitem{CW} R. R. Coifman and G. Weiss, 
 {\it Transference methods in analysis}, 
 CBMS Regional Conference Series in Mathematics 31, 
 Amer. Math. Soc., Providence, 1977.  


\bibitem{CGZ} T. Coulhon, A. Grigor'yan and F. Zucca, 
 `The discrete integral maximum principle and its applications', 
 Tohoku Math. J. 57 (2005), 559-587.      




\bibitem{CSal} T. Coulhon and L. Saloff-Coste, 
 `Puissances d'un op\'{e}rateur r\'{e}gularisant', 
Ann. Inst. H. Poincar\'{e} Probab. Statist. 26 (1990), 
419-436.  



\bibitem{Dav2} E. B. Davies, 
 {\it Heat kernels and spectral theory}, 
 Cambridge Tracts in Math. 92, Cambridge Univ. Press, Cambridge, 1989.   

\bibitem{Dav10} E. B. Davies, 
 `Uniformly elliptic operators with measurable coefficients', 
 J. Funct. Anal. 132 (1995), 
 141-169.  

\bibitem{Dun10} N. Dungey, 
 `On Gaussian kernel estimates on groups', 
 Coll. Math. 100 (2004), 77-90.   

\bibitem{Dun15} N. Dungey, 
 `Heat kernel and semigroup estimates for sublaplacians with 
 drift on Lie groups', 
 Publ. Mat. 49 (2005), 375-391.  

\bibitem{Dun20} N. Dungey, 
 `Properties of random walks on discrete groups: time regularity 
 and off-diagonal estimates',  
 preprint, University of New South Wales, 2005.   

\bibitem{Dun30} N. Dungey, 
 `A note on time regularity for discrete time heat kernels', 
Semigroup Forum 72 (2006), 
404-410. 


\bibitem{Grig} A. Grigor'yan, 
 `Estimates of heat kernels on Riemannian manifolds', 
 in: 
 {\it Spectral theory and Geometry. ICMS Instructional Conference, 
 Edinburgh 1998}, eds. E. B. Davies and Y. Safarov, 
 London Math. Soc. Lecture Note Series 273, 
 Cambridge Univ. Press, 1999, 140-225. 

\bibitem{Grigc} R. Grigorchuk, 
 `Degrees of growth of finitely generated groups and the theory 
of invariant means', 
 Izv. Akad. Nauk. SSSR Ser. Mat. 48 (1984), 
 939-985. 


\bibitem{Gui1} Y. Guivarc'h, 
 `Croissance polynomiale et p\'{e}riodes des fonctions 
  harmoniques', 
 Bull. Soc. Math. France 101 (1973), 333-379.   

\bibitem{Heb1} W. Hebisch, 
 `On heat kernels on Lie groups', 
 Math. Z. 210 (1992), 593-605.    


\bibitem{HebSal} W. Hebisch and L. Saloff-Coste, 
 `Gaussian estimates for Markov chains and random walks on groups', 
 Ann. Probab. 21 (1993), 
 673-709.   

\bibitem{HR} E. Hewitt and K. A. Ross, 
 {\it Abstract Harmonic Analysis I}, 
 Second edition, Grundlehren Math. Wiss. 115, 
 Springer-Verlag, New York, 1979.  


\bibitem{Kat1} T. Kato,  
 {\it Perturbation theory for linear operators},  
 Second edition, Grundlehren Math. Wiss. 132, 
 Springer-Verlag, Berlin, 1984.   


\bibitem{Mac1} G. Mackey, 
 `Induced representations of locally compact groups I', 
 Ann. of Math. 55 (1952), 101-139.  
 

\bibitem{Mel} C. Melzi, 
 `Large time estimates for heat kernels in nilpotent Lie groups', 
 Bull. Sci. Math. 126 (2002), 71-86.  


\bibitem{MZ} D. Montgomery and L. Zippin, 
 {\it Topological transformation groups}, 
 Interscience Publishers, New York-London, 1955. 



\bibitem{Mus} S. Mustapha, 
 `Gaussian estimates for heat kernels on Lie groups', 
 Math. Proc. Cambridge Philos. Soc. 128 (2000), 
 45-64. 



\bibitem{N1} O. Nevanlinna, 
 {\it Convergence of iterations for linear equations}, 
 Birkh\"{a}user, Basel, 1993. 

\bibitem{N2} O. Nevanlinna, 
 `On the growth of the resolvent operators for power bounded operators', 
 in: 
 {\it Linear Operators}, J. Janas, F. H. Szafraniec and J. Zemanek 
(eds.), 
Banach Center Publ. 38, 
Inst. Math., Polish Acad. Sci., 
1997, 247-264.   
v

\bibitem{Pat} A. Paterson, 
 {\it Amenability}, 
 Mathematical Surveys and Monographs 29, 
 American Math. Soc., 
 Providence, 1988.  

 

\bibitem{Robm} D. W. Robinson, 
 {\it Elliptic operators and Lie groups}, 
 Oxford Mathematical Monographs, 
 Oxford University Press, 
 Oxford, 1991.   



\bibitem{Vara} V. S. Varadarajan, 
{\it Lie groups, Lie algebras and their representations}, 
 Graduate Texts in Mathematics 102, 
Springer-Verlag, New York etc., 1984.   



\bibitem{Var10} N. T. Varopoulos, 
 `Long range estimates for Markov chains', 
 Bull. Sci. Math. 109 (1985), 
 225-252.   



\bibitem{Var6} N. T. Varopoulos, 
 `Analysis on Lie groups', 
Rev. Mat. Iberoamericana 12 (1996), 791-917.   







\bibitem{Var20} N. T. Varopoulos, 
 `Distance distortion on Lie groups', 
 in: 
 {\it Random walks and discrete potential theory (Cortona, 1997)},  
 M. Picardello and W.Woess (eds.), 
 Sympos. Math. XXXIX, 
 Cambridge Univ. Press, Cambridge, 1999, 
 320-357.     


\bibitem{Var30} N. T. Varopoulos, 
  `Geometric and potential theoretic results on Lie groups', 
 Canad. J. Math. 52 (2000), 412-437.   


\bibitem{VSC} N. T. Varopoulos, L. Saloff-Coste and T. Coulhon,  
 {\it Analysis and geometry on groups}, 
 Cambridge Tracts in Mathematics 100, 
 Cambridge University Press, Cambridge, 1992.  


\end{thebibliography}
\end{document}